\documentclass{amsart}
\input epsf

\usepackage{longtable}
\newtheorem{thm}{Theorem}[section]

\newtheorem{prop}[thm]{Proposition}
\theoremstyle{definition}
\newtheorem{defn}[thm]{Definition}
\theoremstyle{remark}
\newtheorem{rem}[thm]{Remark}
\newtheorem{rems}[thm]{Remarks}
\numberwithin{equation}{section}
\newcommand{\R}{\mathbb R}

\newcommand{\C}{\mathbb{C}}
\newcommand{\Z}{\mathbb{Z}}
\newcommand{\inv}{^{-1}}
\newcommand{\si}{\sigma}
\begin{document}

\title[Intersection of opposed real big cells for $G_2$]
{The intersection of opposed big cells in the real flag variety
of type $G_2$}
\author{R. J. Marsh}%
\address{University of Leicester}%
\email{R.Marsh@mcs.le.ac.uk}%

\author{K. Rietsch}%
\address{DPMMS, Cambridge University and Newnham College}%
\email{rietsch@dpmms.cam.ac.uk}%

\thanks{The second named author was supported by EPSRC
fellowship GR/M09506/01 and both authors were supported by a
University of Leicester Research Fund
Grant.}%
\subjclass{}%
\keywords{}%

\date{}%
%\dedicatory}%
%\commby}%
% ----------------------------------------------------------------
\begin{abstract}
We compute the Euler characteristics of the individual  connected
components of the intersection of two opposed big cells in the
real flag variety of type $G_2$, verifying a conjecture from
\cite{Rie:BruCel}.
%There are $11$ connected components, one of which
%turns out to have Euler characteristic $2$ while the others have
%Euler characteristic $1$.
\end{abstract}
\maketitle
% ----------------------------------------------------------------
 \section{Introduction}

Let $G$ be a connected linear algebraic group.  Any flag variety
(homogeneous projective variety $G/B$) has a myriad of cell
decompositions, so--called `Bruhat decompositions'. For every
Borel subgroup there is precisely one, and the cells in the
decomposition are simply the orbits under this group. Furthermore
there is always a unique open dense orbit called the `big cell'.
In the present paper we fix two Borel subgroups opposite to one
another and study the intersection of the two resulting big
cells. More precisely we are interested in the real points of
this variety, in the case where everything is split over $\R$.

We let $\mathcal B^*$ denote the real points of the intersection
of two opposed big cells, endowed with the Hausdorff topology
coming from $\R$. There has been some recent work on determining
the number of connected components
%(parameterizing them by the
%connected components of a combinatorially defined graph)
and the
Euler characteristics
of these varieties $\mathcal B^*$, see e.g. \cite{Rie:BruCel,
Rie:BigCel}, and for the type $A$ case \cite{SSV:BigCel1,
SSV:BigCel2}. In particular it was open for a long time whether
the connected components are always contractible. This has proved
not to be the case. The first example is in type $G_2$ found in
\cite{Rie:BruCel} (but it is also false in type $A$ for rank $>4$
as follows from the explicit formula for the number of connected
components in \cite{SSV:BigCel2} and a computation of the Euler
characteristics e.g. using \cite{Rie:BruCel}).

Therefore to the graph defined in \cite{Rie:BruCel}
parameterizing the connected components of
$\mathcal B^*$ one can in principle add
another nontrivial datum. That is, every component of the graph
has an integer associated to it given simply by the Euler
characteristic of the corresponding connected component of
$\mathcal B^*$. Lusztig showed in
\cite{Lus:BruCelApp} that the graphs arising in this setting are
`mod $2$' quotients of the
graphs parameterizing the canonical basis of the corresponding
quantum enveloping algebra $\mathcal U^-$. Therefore by taking the
Euler characteristics one is assigning in a natural way an integer
to every canonical basis element. No direct canonical basis
interpretation of these integers is known. Moreover they
have not up to now been computed in any nontrivial examples.

We focus our attention here on the $G_2$ case. In Figure~\ref{f:oldfig}
we have reproduced the parameterization
\cite{Rie:BruCel} of the connected components
of $\mathcal B^*$ in type $G_2$. There are $11$
connected components, and one is
conjectured to have Euler characteristic $2$ (the last one),
while the others should have Euler characteristic $1$. By an
independent method the total Euler characteristic of $\mathcal
B^*$ was worked out in \cite{Rie:BruCel} as $12$.

The aim of this paper is to compute the individual Euler
characteristics of the connected components of $\mathcal B^*$ to
verify the above conjecture. We use two main tools. The first one
is a cell decomposition of $\mathcal B^*$ due to Deodhar
\cite{Deo:Decomp}, which is recalled in the next section. And the
second ingredient is Berenstein and Zelevinsky's Chamber Ansatz,
explained in Section~3, which we require to be able to tell which
connected component each of the various cells of Deodhar's
decomposition lie in. The resulting decompositions for the
connected components are indicated in Figure~\ref{components}.
The Euler characteristics are the expected ones as can be read off
from that figure.

\section{Decomposition of $\mathcal B^*$}

\subsection{Preliminaries}

We begin by introducing some notation and standard facts as can
be found e.g. in \cite{Springer:book}. Let $G_\C$ be a complex
simple linear algebraic group. After this general section we will
choose it specifically as the one of type $G_2$. In any case we
always focus on its split real form $G_\R=G$. All the varieties
in this paper will be defined over $\R$, and we identify them
with their $\R$--valued points. The topology we consider is the
usual Hausdorff topology coming from $\R$ (rather than Zariski
topology). We write $\R^*$ for $\R\setminus\{ 0\}$.

Let $\mathcal B$ denote the (real) flag variety of $G$. The
elements of $\mathcal B$ are the Borel subgroups of $G$. We
sometimes write $[B]$ for the Borel subgroup $B$ considered as a
point in the flag variety. The transitive action of $G$ on
$\mathcal B$ is denoted by $g\cdot [B]:=[g B g\inv]$.

Let $T\subset G$ be a fixed split maximal torus, and $B^+$ a fixed
Borel subgroup which contains $T$. We also automatically have
given the opposite Borel subgroup $B^-$ (such that $B^+\cap
B^-=T$) and the unipotent radicals $U^+$ and $U^-$ of these two
Borel subgroups.

Let $X^*(T)$ and $X_*(T)$ be the character, respectively
cocharacter lattices of $T$ with their canonical perfect pairing
\begin{equation*}
\left<\ ,\ \right >:\ X_*(T)\times X^*(T)\to\Z.
\end{equation*}
And let $A=(A_{ij})$ be the Cartan matrix of $G$. We denote the
positive simple roots (corresponding to $B^+$) by
$\alpha_1,\dotsc, \alpha_r\in X^*(T)$. We will write $\alpha >0$
if $\alpha\in X^*(T)$ is a positive root, i.e. a nonnegative
linear combination of the simple roots. We also have the simple
coroots $\alpha_1^\vee,\dotsc, \alpha_r^\vee\in X_*(T)$ which are
determined by $\left <\alpha_i^\vee,\alpha_j\right > = A_{ij}$.

We will make extensive use of the $1$--parameter simple root
subgroups: Let us fix for $i=1,\dotsc, r$, Chevalley generators
$e_i$ and $f_i\in Lie(G)$ of the (real) Lie algebra, with the
$e_i$'s lying in positive simple root spaces. Then we define
\begin{eqnarray*}
x_i:\R\to B^+,&\qquad x_i(t):=\exp(t e_i),\\
y_i:\R\to B^-,&\qquad y_i(t):=\exp(t f_i).
\end{eqnarray*}

The Weyl group $N_G(T)/T$ of $G$ is denoted $W$ and is generated
by the usual simple reflections $s_1,\dotsc, s_r$ ($r$ being the
rank of $G$). We fix a representative $\dot s_i\in G$ for the
Weyl group element $s_i$ by defining $\dot s_i:=x_i(1) y_i(-1)
x_i(1)$. For any other $w\in W$ choose a reduced (minimal number
of factors) expression $w=s_{i_1}\dotsc s_{i_k}$ and set $\dot
w:=\dot s_{i_1}\dotsc \dot s_{i_k}$ to get a representative. It
is well known that this definition of $\dot w$ is independent of
the reduced expression.
%(see e.g. \cite{Springer:book}).
Note that
$x_i(t)\in  U^+\cap B^-\dot s_i B^-$ and $y_i(t)\in  U^-\cap
B^+\dot s_i B^+$, whenever $t\ne 0$. Explicitly,
\begin{eqnarray*}
x_i(t)&=&y_i(t\inv)\ \dot s_i\alpha_i^\vee(-t)\ y_i(t\inv ),\\
y_i(t)&=&x_i(t\inv)\ \alpha_i^\vee(-t)\dot s_i \ x_i(t\inv ).
\end{eqnarray*}

The Weyl group acts on $T$ and hence naturally also on $X_*(T)$
and $X^*(T)$. The length of (a reduced expression of) an element
$w\in W$ is denoted $\ell(w)$. We also consider the Bruhat order
$<$ on $W$. The statement $ws_i<w$ is equivalent to
$w\cdot\alpha_i<0$, that $w$ sends $\alpha_i$ to a negative root.

There is a unique longest element of
$W$ which is called $w_0$. We propose to study the open
subvariety of the flag variety
\begin{equation*}
\mathcal B^*:=B^+\dot w_0\cdot [B^+]\cap B^-\dot w_0\cdot [B^-] \
\subset \ \mathcal B,
\end{equation*}
in other words the intersection of the two big cells for the
Bruhat decompositions relative to $B^-$, respectively $B^+$.

Note that $\mathcal B^*$ is an affine variety and can be
identified naturally with open subvarieties both of $U^+$ and of
$U^-$. We have two isomorphisms
\begin{equation}\label{e:isos}
\begin{matrix}
&&\mathcal B^* &&    \\
&i_+\nearrow & &\nwarrow i_-&                                   \\
U^+\cap(B^-\dot w_0 B^-)&&&&U^-\cap(B^+\dot w_0 B^+),
\end{matrix}
\end{equation}
where $i_+$ takes $u\in U^+$ to $u\cdot[B^-]$ while the right hand
map $i_-$ takes $u\in U^-$ to $u\cdot[B^+]$.

\subsection{Deodhar's decomposition of $\mathcal B^*$}

We begin by making the following definition (see
\cite[Appendix]{KaLus:Hecke}).

\begin{defn}(relative position in $\mathcal B$)
Let $w\in W$ and consider the action of $G$ on $\mathcal
B\times\mathcal B$ by simultaneous conjugation. We say that two
Borel subgroups $B, B'\in \mathcal B$ are in relative position
$w$ if the pair $(B,B')$ lies in the $G$--orbit of $(B^+,\dot
w\cdot B^+)$. From Bruhat decomposition it follows that such $w$
is unique and exists for any pair $(B, B')$. We write $B\overset
w\longrightarrow B'$.
\end{defn}

If $B\overset w\longrightarrow B'$ and $w=s_{i_1}\dotsc s_{i_n}$
is a reduced expression, then it also follows from Bruhat
decomposition that there exist uniquely determined Borel subgroups
$B_0, B_1, \dotsc, B_n=B'$ such that
\begin{equation*}
B=B_0\overset {s_{i_1}}\longrightarrow
B_1\overset{s_{i_2}}\longrightarrow B_2\longrightarrow
\cdots\longrightarrow B_{n-1} \overset{s_{i_n}}\longrightarrow
B_n=B'
\end{equation*}
Explicitly, $B\overset w\longrightarrow B'$ means $B=g\cdot B^+$
and $B'=g\dot w\cdot B^+$ for some $g\in G$. And the $B_j$'s are
then given by $B_j=g\dot s_{i_1}\dotsc\dot s_{i_j}\cdot B^+$.

Using the above notation, we have $\mathcal B^*=\{B\in \mathcal B\
|\ B^+\overset{w_0}\longrightarrow B\overset{w_0}\longleftarrow
B^-\}$. The following definition is a special case of
\cite[Definition 2.3]{Deo:Decomp}.

\begin{defn}[distinguished subexpressions for $1$]
Let $s_{i_1}\dotsc s_{i_n}=w\in W$ be a reduced expression for
$w$ (so $n=\ell(w)$). Then by a {\sl subexpression} of this
reduced expression we mean a sequence of Weyl group elements
$\si=(\si_0,\si_1,\dotsc, \si_n)$ such that
\begin{equation*}
\si_0=1, \qquad \si_j=\begin{cases} \text{either }&\si_{j-1}s_{i_j} \\
\quad \text{or}&\si_{j-1},
\end{cases}\qquad \text{for all $j=1,\dotsc n$.}
\end{equation*}
In particular the ``empty'' subexpression $\si=(1,\dotsc, 1)$ is
allowed. We call $\si$ a {\sl subexpression for $1$} if $\si_n=1$.

A subexpression is called {\sl distinguished} if  we have
\begin{equation*}
\si_j\le \si_{j-1}\ s_{i_j},\quad \text{for all
$j\in\{1,\dotsc,n\}$}.
\end{equation*}
This means that if right  multiplication by $s_{i_j}$ decreases
the length of $\si_{j-1}$, then we must choose
$\si_{j}=\si_{j-1}s_{i_j}$ to get a distinguished subexpression.
\end{defn}

The following result is stated by Deodhar over an algebraically
closed field, but it extends trivially to any split form, so we
state it here in our present setting over the reals.

\begin{thm}[\cite{Deo:Decomp}]\label{t:Deo}
Let $s_{i_1}\dotsc s_{i_N}$ be a fixed reduced expression for
$w_0$.
\begin{enumerate}
\item
Suppose $B=x\dot w_0\cdot B^+$ is an element of $\mathcal B^*$,
where $x\in B^+$. Then the sequence $(\si_0,\dotsc,\si_N)=:\si(B)$
defined by
\begin{equation*}
x s_{i_1}\dotsc s_{i_k}\in B^-\si_{k} B^+
\end{equation*}
is a well--defined distinguished subexpression for $1$ (of
$s_{i_1}\dotsc s_{i_N}$).
\item
Let $D_\sigma:=\{ B\in\mathcal B^*\ |\ \si(B)=\si \}$. Then
\begin{equation*}
 \begin{array}{lcl} D_\sigma\cong (\R^*)^{|I(\sigma)|}\times \R^{|K(\sigma)|},
 &
  \text{ where }& I(\sigma)=\{j\in\{1,\dotsc, N\}\ |\
\si_j=\si_{j-1}\},\\
&\text{ and }& K(\sigma)=\{j\in\{1,\dotsc, N\}\ |\
\si_j<\si_{j-1}\}.
\end{array}
\end{equation*}
\end{enumerate}
\end{thm}

An isomorphism as in part (2) of the theorem (albeit not identical to
the one in \cite{Deo:Decomp}) will be constructed
explicitly below. We continue to fix the reduced expression
$s_{i_1}\dotsc s_{i_N}$ for $w_0$ in what follows.
\begin{rem}
The definition of $D_\sigma$ is natural to state using relative
position. If $B\in\mathcal B^*$, and $B_1,\dotsc, B_N$ are defined
by
\begin{equation*}
B^+\overset {s_{i_1}}\longrightarrow
B_1\overset{s_{i_2}}\longrightarrow B_2\longrightarrow
\cdots\longrightarrow B_{N-1} \overset{s_{i_N}}\longrightarrow
B_N=B,
\end{equation*}
then $B\in D_\sigma$ precisely if $B^-\overset{
w_0\sigma_j}\longrightarrow B_j$ for all $j$. So the $\si$
controls the relative positions of the intermediate $B_i$'s with
respect to $B^-$.
\end{rem}

\subsection{Inductive construction}\label{s:induction}
Let $\si$ be a fixed distinguished subexpression for $1$ of
$s_{i_1}\dotsc s_{i_N}$. We now want to describe
$D_{\sigma}^-\subset U^-\cap B^+ w_0 B^+$, the preimage of
$D_{\sigma}$ under the isomorphism $i_-$ in \eqref{e:isos}. We
will do this by building up the elements $B\in D_\si$ from the
$B_j$'s defined in the previous remark. To begin with note that
if $B_{j-1}=g\cdot B^+$, then
$B_{j-1}\overset{s_{i_j}}\longrightarrow B_j$ just says that
$B_{j}=g x_{i_j}(t)\dot s_{i_j}\cdot B^+$ for some $t\in\R$, or
equivalently
\begin{equation*}
B_{j}=\begin{cases}\text{either} & g y_{i_j}(t)\cdot B^+\quad \text{some $t\in\R^*$,}\\
\text{ or} &g\dot s_{i_j}\cdot B^+.
\end{cases}
\end{equation*}

Let us first determine the possible $B_1$'s. There are two cases.
\begin{itemize}
\item
If $\si_1=1$, then we have $B^+\overset{s_{i_1}}\longrightarrow
B_1\overset{w_0}\longleftarrow B^-$. Therefore
$B_1=y_{i_1}(t)\cdot B^+$ for some $t\in\R^*$.
\item
If on the other hand $\si_1=s_{i_1}$ then
$B^+\overset{s_{i_1}}\longrightarrow B_1\overset{w_0
s_{i_1}}\longleftarrow B^-$ and we get $B_1=\dot s_{i_1}\cdot
B^+$.
\end{itemize}

Suppose in general we have $B_{j-1}=g\cdot B^+$ given, where
$g=y\dot \si_{j-1}$ for some $y\in U^-$.
%Now let us assume by induction that $B_{j-1}=g\cdot B^+$ where
%$g= z_1\dotsc z_{j-1}$ and  $z_k$ is either of the form
%$y_{i_k}(t)$ for some $t\in\R^*$ (if $\si_{k-1}=\si_k$), or
%$z_k=\dot s_{i_k}$ (if $\si_{k-1}<\si_k$), or $z_k=x_{i_k}(m)\dot
%s_{i_k}$ some $m\in\R$ (if $\si_{k-1}>\si_k$).
We then want to construct all possible $B_j=g'\cdot B^+$ from
this $B_{j-1}$. There are three cases.
\begin{enumerate}
\item
Suppose first that $\si_j=\si_{j-1}$. Then we have the setting
\begin{equation*}
\begin{matrix}
 B^- &\overset{w_0\si_{j}}\longrightarrow &B_{j-1}\\
    {}_{w_0\si_{j}}\searrow& &\swarrow_{s_{i_j}}\\
                              &  B_j&
\end{matrix}
\end{equation*}
So if $B_{j-1}=g\cdot B^+$ where $g=y\dot \si_{j-1}$, then it is
easy to rule out $g\dot s_{i_j}\cdot B^+$ for $B_j$ and the only
possible solutions are of the form $B_j=gy_{i_j}(t)\cdot B^+$ for
some $t\in\R^*$.
%And this exhausts all possible $B_j$'s since
%the only other candidate, $g\dot s_{i_j}\cdot B^+$, has the wrong
%relative position with respect to $B^-$.
Note that $w_0\si_{j} s_{i_j} < w_0\si_{j} $ (since $\si$ is
distinguished). We claim that therefore $B_j=gy_{i_j}(t)\cdot
B^+$ has the correct relative positions for any $t\in\R^*$ and in
fact $g':=gy_{i_j}(t)\in U^- \dot\si_j$. All of this follows
since we have $\si_{j-1}\cdot\alpha_{i_j}>0$ and therefore
\begin{equation*}
gy_{i_j}(t)=y\dot \si_{j-1} y_{i_j}(t)=y\dot
\si_{j-1}y_{i_j}(t)\dot\si_{j-1}\inv\dot\si_j\in U^-\dot\si_j,
\end{equation*}
using $\dot\si_{j-1}=\dot\si_{j}$.
\item
Suppose next that $\si_j>\si_{j-1}$. Then
\begin{equation*}
\begin{matrix}
 B^- &\overset{w_0\si_{j-1}}\longrightarrow &B_{j-1}\\
    {}_{w_0\si_{j}}\searrow& &\nearrow_{s_{i_j}}\\
                              &  B_j&
\end{matrix}
\end{equation*}
and since the lengths add,
$\ell(w_0\si_{j})+\ell(s_{i_j})=\ell(w_0\si_{j-1})$, we get that
$B_j$ is uniquely determined by $B_{j-1}$ and equals to $g\dot
s_{i_j}\cdot B^+$. We immediately have that $g':=g\dot s_{i_j}\in
U^-\dot\si_j$.
\item
Finally we have the case $\si_j<\si_{j-1}$. In this case the
other two lengths add, $\ell(w_0\si_{j-1})+\ell(s_{i_j})=
\ell(w_0\si_j)$, and the diagram
\begin{equation*}
\begin{matrix}
 B^- &\overset{w_0\si_{j-1}}\longrightarrow &B_{j-1}\\
    {}_{w_0\si_{j}}\searrow& &\swarrow_{s_{i_j}}\\
                              &  B_j&
\end{matrix}
\end{equation*}
is automatically satisfied for any $B_j$ in position $s_{i_j}$
relative to $B_{j-1}$. Therefore we can take $B_j=g x_{i_j}(m)\dot
s_{i_j}\cdot B^+$ for any $m\in \R$.

We claim that in this case $g':=g x_{i_j}(m)\dot s_{i_j}\inv \in
U^-\dot \si_j$. This holds because we have $\si_{j-1}\cdot
\alpha_{i_j}<0$. Therefore $\dot \si_{j-1}x_{i_j}(t)\dot
\si_{j-1}\inv\in U^-$ and
\begin{equation*}
gx_{i_j}(m)\dot s_{i_j}\inv=y\dot \si_{j-1} x_{i_j}(m)\dot
s_{i_j}\inv=y\dot \si_{j-1}x_{i_j}(m)\dot \si_{j-1}\inv \dot
\si_{j}\in U^-\dot\si_j.
\end{equation*}
\end{enumerate}

Applying this procedure recursively to express finally $B_N$ we
get the following proposition (and we also recover Deodhar's
result Theorem~\ref{t:Deo}).

\begin{prop}\label{p:explicit} For $w_0=s_{i_1}\dotsc s_{i_N}$ a fixed reduced
expression, and $\si$ a distinguished subexpression for $1$ let
\begin{eqnarray*}
 I(\sigma)=\{j\,:\,1\leq j\leq N\mbox{\ and\
 }\sigma_{j-1}=\sigma_j\},\\
 J(\sigma)=\{j\,:\,1\leq j\leq N\mbox{\ and\
 }\sigma_{j-1}<\sigma_j\},\\
 K(\sigma)=\{j\,:\,1\leq j\leq N\mbox{\ and\
 }\sigma_{j-1}>\sigma_j\}.
\end{eqnarray*}
Then we have explicitly
\begin{equation*}
D_\sigma^-=\left\{z_1 z_2\dotsc z_N\ \left | \
z_j=\begin{cases}y_{i_j}(t_j)&\text{ if $j\in I(\sigma)$}\\
\dot s_{i_j}&\text{ if $j\in J(\sigma)$}\\
x_{i_j}(m_j)\dot s_{i_j}\inv&\text{ if $j\in K(\sigma)$}
\end{cases}\text{ , where $m_j\in\R$, $t_j\in\R^*$}\right . \right\}.
\end{equation*}
\end{prop}

\begin{proof}
By the inductive construction in \ref{s:induction} above we showed
that any $B\in D_{\si}$ is of the form $z\cdot B^+$ for $z=z_1
\dotsc z_N$ as in the proposition (and that these $z\cdot B$'s
all lie in $D_{\si}$). We also showed that $z_1\dotsc z_k\in
U^-\dot\si_{k}$ for all $k$, so in particular $z_1\dotsc z_N\in
U^-$, since $\si_N=1$. Therefore $D^-_{\sigma}$ is precisely the
set of these $z=z_1\dotsc z_N$.
\end{proof}

\begin{rem} Note that the map
\begin{equation*}
(\R^*)^{I(\si)}\times \R^{K(\si)}\to D^-_{\sigma}
\end{equation*}
arising from the proposition is an isomorphism, with inverse
basically constructed in \ref{s:induction}.
\end{rem}

\subsection{Further refinement}\label{s:cells}

Over $\R$ it is natural to consider the connected components of
the $D_\sigma$'s to get a cell decomposition. Let
$h:I(\si)\rightarrow \{1,-1\}$ be a choice of signs for the
elements of $I(\sigma)$. Then we define
\begin{equation*}
 D_{\sigma}^-(h)=\left\{z_1
z_2\dotsc z_N\
\left |\ z_j=\begin{cases}y_{i_j}(t_j)\text{ with $h(j)t_j\in\R_{>0}$}&\text{ if $j\in I(\sigma)$,}\\
\dot s_{i_j}&\text{ if $j\in J(\sigma)$,}\\
x_{i_j}(m_j)\dot s_{i_j}^{-1}\text{ with  $m_j\in\R $}&\text{ if $j\in
K(\sigma)$.}
\end{cases}\right . \right\}.
%\{\prod_{j=1}^N p_j(a_j)\,:\,f(j)a_j\in \mathbb{R}_{>0},
%\,j=1\ldots N\}.
\end{equation*}
$D^-_\sigma(h)$ is a (real) semi--algebraic cell in $U^-\cap B^+
\dot w_0 B^+$ of dimension $|I(\si)|+|K(\si)|$. Its image
$i_-(D^-_\sigma(h))$ in $\mathcal B^*$ is denoted $D_\sigma(h)$.

\section{Type $G_2$}

\begin{figure}\label{f:roots}
\epsfbox{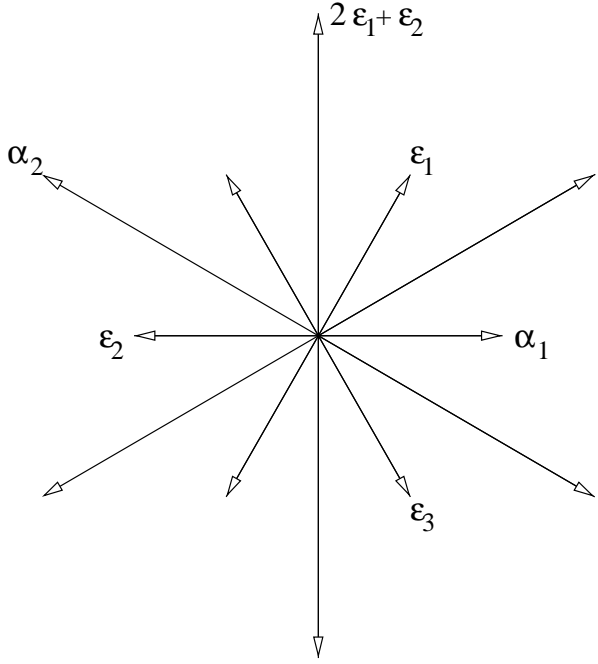} \caption{$G_2$ root system}
\end{figure}

From now on let $G$ be of type $G_2$. Then the Cartan matrix
$A=(A_{ij})$ for type $G_2$ is given by
\begin{equation*}
A=\begin{pmatrix} 2 &-3\\
-1 &2
\end{pmatrix}.
\end{equation*}
The Weyl group $W$ has two generators $s_1, s_2$ corresponding to
reflection by the short root $\alpha_1$ and the long root
$\alpha_2$, respectively (see Figure \ref{f:roots}). The
fundamental weights are $\omega_1=\varepsilon_1$, giving rise to
a $7$--dimensional representation, and $\omega_2=2\varepsilon_1
+\varepsilon_2$, the highest weight of the $14$--dimensional
adjoint representation.

The longest element is
\begin{equation*}
w_0=s_1 s_2 s_1 s_2 s_1 s_2=s_2 s_1 s_2 s_1 s_2 s_1.
\end{equation*}
We let $\mathbf i=(1,2,1,2,1,2)$ and $\mathbf{\widetilde
i}=(2,1,2,1,2,1)$ stand for these two reduced expressions of
$w_0$.

To give Deodhar's decomposition of $\mathcal B^*$ in this case we
first need to fix a reduced expression of $w_0$, so let us pick
$\mathbf i=(1,2,1,2,1,2)$. Next we list all of the distinguished
subexpressions for $1$ of $\mathbf i$:
\begin{equation*}
\Sigma_{\mathbf i}:=\left \{
\begin{array}{l  l}
 xxxxxx :=(1,1,1,1,1,1,1),& xxx2x2:=(1,1,1,1,s_2,s_2,1),\\
 1x1xxx :=(1,s_1,s_1,1,1,1,1),& 1x12x2:=(1,s_1,s_1,1,s_2,s_2,1),\\
 x2x2xx :=(1,1,s_2,s_2,1,1,1),&12x21x:=(1,s_1,s_1 s_2,s_1 s_2, s_1,
 1,1),\\
  xx1x1x :=(1,1,1,s_1,s_1,1,1),
 & x21x12:=(1,1,s_2,s_2 s_1, s_2 s_1,
 s_2,1)
\end{array}\right \}.
\end{equation*}
Furthermore, define
\begin{eqnarray*}
y_{xxxxxx}(t_1,\dotsc,t_6):=y_{\mathbf i}
(t_1,\dotsc,t_6)&:=&y_1(t_1)y_2(t_2)y_1(t_3)
y_2(t_4)y_1(t_5)y_2(t_6)\\
y_{1x1xxx}(t_1,m_1,t_2,t_3,t_4)&:=& \dot s_1 y_2(t_1) x_1(m_1)
\dot s_1\inv
y_2(t_2) y_1(t_3) y_2(t_4)\\
y_{x2x2xx}(t_1,t_2,m_1,t_3,t_4)&:=&y_1(t_1) \dot s_2 y_1(t_2)
x_2(m_1)\dot s_2\inv
y_1(t_3) y_2(t_4)\\
y_{xx1x1x}(t_1,t_2,t_3,m_1,t_4)&:=&y_1(t_1) y_2(t_2) \dot s_1
y_2(t_3) x_1(m_1)\dot s_1\inv y_2(t_4)\\
y_{xxx2x2}(t_1,t_2,t_3,t_4,m_1)&:= &y_1(t) y_2(t_2) y_1(t_3)\dot
s_2
y_1(t_4) x_2(m_1) \dot s_2\inv\\
y_{1x12x2}(t_1,m_1,t_2,m_2)&:=&\dot s_1 y_2(t_1) x_1(m_1)\dot
s_1\inv \dot s_2
y_1(t_2) x_2(m_2)\dot s_2\inv\\
y_{12x21x}(t_1,m_1,m_2,t_2)&:=&\dot s_1\dot s_2 y_1(t_1) x_2(m_1)
\dot s_2\inv
x_1(m_2) \dot s_1\inv y_2(t_2)\\
y_{x21x12}(t_1,t_2,m_1,m_2)&:=&y_1(t_1)\dot s_2\dot s_1 y_2(t_2)
x_1(m_1) \dot s_1\inv x_2(m_2) \dot s_2\inv.
\end{eqnarray*}
By Proposition \ref{p:explicit} we have
\begin{eqnarray*}
% D_{\pm,\pm ,\pm ,\pm ,\pm ,\pm}=
D_{\mathbf i}^-:=D_{xxxxxx}^-&=&\{y_{\mathbf i}(t_1,t_2,t_3,t_4,t_5,t_6)\ |\ t_1,\dotsc, t_6\in\R^*\}\\
% D_{0,\pm, *, \pm, \pm, \pm}=
D_{1x1xxx}^-&=&\{y_{1x1xxx}(t_1,m_1,t_2,t_3,t_4)\ |\ t_1,\dotsc,
    t_4 \in\R^*,\  m_1\in\R \}\\
&&\vdots \\
%D_{\pm,0,0,\pm,*,*}=
D_{x21x12}^-&=&\{y_{x21x12}(t_1,t_2,m_1,m_2)\ |\ t_1,t_2\in\R^*,
m_1, m_2\in\R\},
%\end{array}
\end{eqnarray*}
with the property that
\begin{equation*}
U^-\cap B^+\dot w_0 B^+=\bigsqcup_{\si\in\Sigma_{\mathbf i}}
D_{\sigma}^-.
\end{equation*}
From these we also get the $D^-_{\si}(h)$'s and $D_{\si}(h)$'s
defined in \ref{s:cells}, where $h:I(\sigma)\to\{1,-1\}$
determines the signs of the $t_i\in\R^*$.

\section{Parameterization of the connected components of $\mathcal B^*$ :
using the opposite reduced expression}\label{s:param}

\begin{center}
\begin{figure}
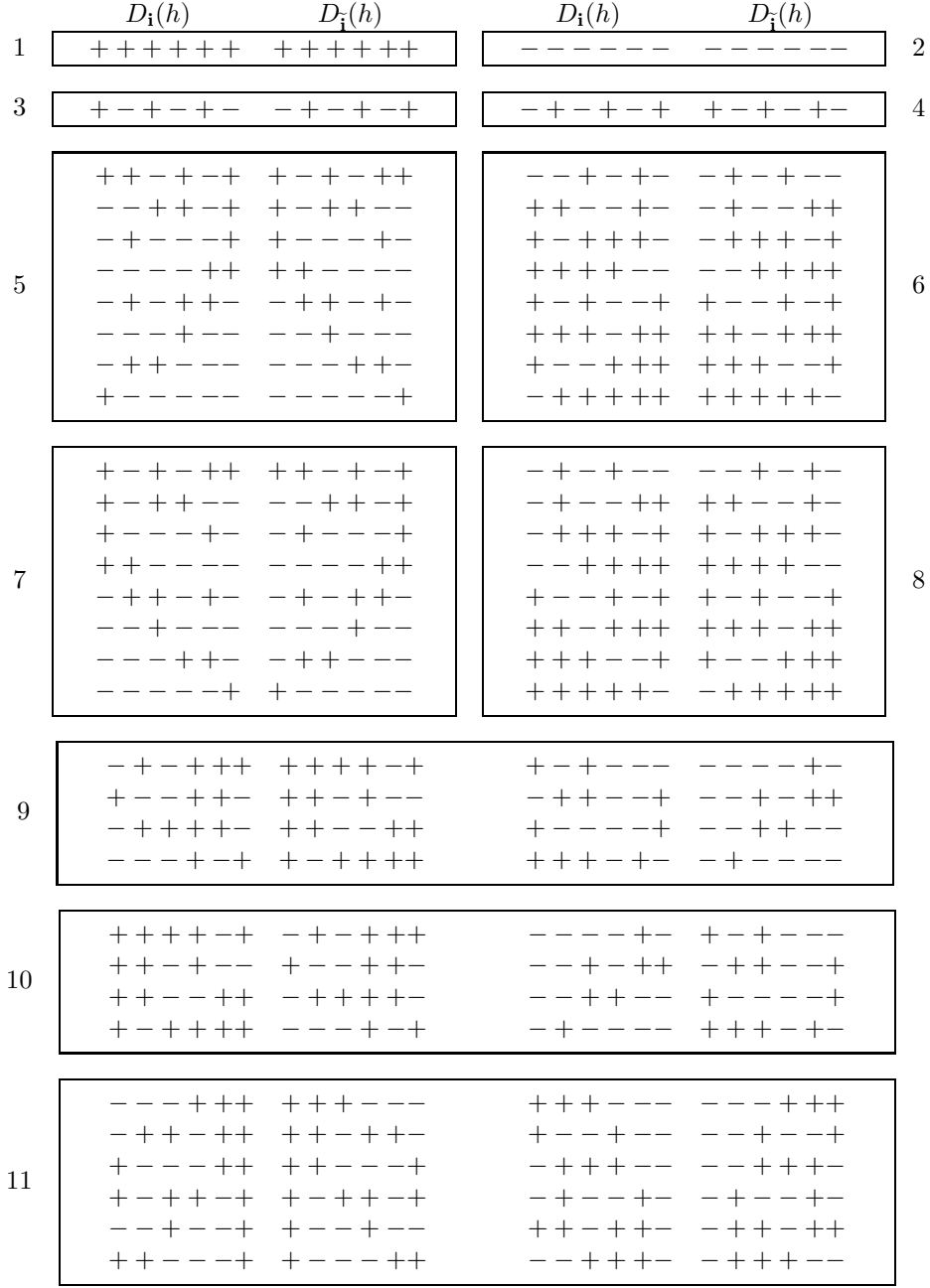

\begin{tabular}{cccc}
 & $D_{\mathbf i}(h)\hskip 1.7cm D_{\widetilde{\mathbf i}}(h)$& $D_{\mathbf
 i}(h)\hskip 1.7cm
 D_{\widetilde{\mathbf i}}(h)$&\\
 \ 1 & \framebox[5.4cm]{$++++++\quad ++++++$} &
 \framebox[5.4cm]{$------\quad ------$}
& 2
\end{tabular}

\vspace*{0.3cm}

\begin{tabular}{cccc}
\ 3 &
 \framebox[5.4cm]{$+-+-+-\quad -+-+-+$}
&
 \framebox[5.4cm]{$-+-+-+\quad+-+-+-$}
& 4
\end{tabular}

\vspace*{0.3cm}

\begin{tabular}{cccc}
\ 5 &
 \framebox[5.4cm]{
$
\begin{array}{cc}
++-+-+ & +-+-++  \\
--++-+ & +-++-- \\
-+---+ & +---+- \\
----++ & ++---- \\
-+-++- & -++-+- \\
---+-- & --+--- \\
-++--- & ---++- \\
+----- & -----+
\end{array}
$ } &
 \framebox[5.4cm]{
$
\begin{array}{cc}
--+-+- &  -+-+-- \\
++--+- &  -+--++\\
+-+++- &  -+++-+\\
++++-- &  --++++\\
+-+--+ &  +--+-+\\
+++-++ &  ++-+++\\
+--+++ &  +++--+\\
-+++++ &  +++++-
\end{array}$ }
& 6
\end{tabular}

\vspace*{0.3cm}

\begin{tabular}{cccc}
\ 7 &
 \framebox[5.4cm]{
$
\begin{array}{cc}
+-+-++ & ++-+-+\\
+-++-- & --++-+\\
+---+- & -+---+\\
++---- & ----++\\
-++-+- & -+-++-\\
--+--- & ---+--\\
---++- & -++---\\
-----+ & +-----
\end{array}$ } &
 \framebox[5.4cm]{
$
\begin{array}{cc}
-+-+-- & --+-+-\\
-+--++ & ++--+-\\
-+++-+ & +-+++-\\
--++++ & ++++--\\
+--+-+ & +-+--+\\
++-+++ & +++-++\\
+++--+ & +--+++\\
+++++- & -+++++\\
\end{array}$ } & 8
\end{tabular}

\vspace*{0.3cm}

\begin{tabular}{ccc}
\ 9 & \framebox[11.2cm]{ $
\begin{array}{cc}
 -+-+++ & ++++-+\\
 +--++- & ++-+--\\
 -++++- & ++--++\\
 ---+-+ & +-++++
\end{array}\hskip 1cm
\begin{array}{cc}
 +-+--- & ----+-\\
 -++--+ & --+-++\\
 +----+ & --++--\\
 +++-+- & -+----
\end{array}
$
} & \  \end{tabular}

\vspace*{0.3cm}

\begin{tabular}{ccc}
10 & \framebox[11.2cm]{ $
\begin{array}{cc}
 ++++-+ & -+-+++\\
 ++-+-- & +--++-\\
 ++--++ & -++++-\\
 +-++++ & ---+-+\\
\end{array}\hskip 1cm
\begin{array}{cc}
 ----+- & +-+---\\
 --+-++ & -++--+\\
 --++-- & +----+\\
 -+---- & +++-+-
\end{array}
$ } & \  \end{tabular}

\vspace*{0.3cm}

\begin{tabular}{ccc}
11 & \framebox[11.2cm]{ $
\begin{array}{cc}
 ---+++ & +++--- \\
 -++-++ & ++-++-\\
 +---++ & ++---+\\
 +-++-+ & +-++-+\\
 --+--+ & +--+--\\
 ++---+ & +---++
\end{array}\hskip 1cm
\begin{array}{cc}
 +++--- &  ---+++\\
 +--+-- &  --+--+\\
 -+++-- &  --+++-\\
 -+--+- &  -+--+-\\
 ++-++- &  -++-++\\
 --+++- &  -+++--
\end{array}$ } & \
\end{tabular}
\caption{\label{f:oldfig} Connected components of $\mathcal B^*$
identified with $U^-\cap B^+\dot w_0 B^+$, following
\cite{Rie:BruCel}}
\end{figure}
\end{center}

We recall the parameterization of the set of connected components
of $\mathcal B^*$ from \cite{Rie:BruCel}. Consider the open subset
$D_{\mathbf i}=\bigsqcup_{h:\{1,\dotsc, 6\}\to \{\pm 1\} }
D_{\mathbf i}(h)$ of $\mathcal B^*$ from above. Explicitly,
\begin{equation*}
D_{\mathbf i}(h)=\left \{y_{\mathbf i}(t_1,\dotsc, t_6)\cdot
[B^+]=y_1(t_1) y_2(t_2)\dotsc y_1(t_5) y_2(t_6)\cdot [B^+]\ |\
 h(i) t_i\in\R_{>0} \right \}.
\end{equation*}
We also have analogs of these for ${\mathbf{\widetilde
i}}=(2,1,2,1,2,1)$,
\begin{eqnarray*}
D_{\widetilde{\mathbf i }}(h)&:=&\left \{y_{\widetilde{\mathbf
i}}(t_1,\dotsc, t_6 )\cdot [B^+]:=y_2(t_1)\dotsc y_1(t_6)\cdot
[B^+]\ |\ h(i) t_i
\in\R_{>0}\right \}\\
D_{\widetilde{\mathbf i}}\quad &:=&\bigsqcup_{h:\{1,\dotsc
,6\}\to\{\pm 1\}} D_{\widetilde{\mathbf i}}(h).
\end{eqnarray*}
It was proved in \cite{Rie:BruCel} that the union
$D^*:=D_{\mathbf i}\cup D_{\widetilde{\mathbf i}}$ has complement
of codimension $\ge 2$ in $\mathcal B^*$. Hence the connected
components of $\mathcal B^*$ correspond bijectively to the
connected components of $D^*$. And these were determined roughly
by checking which of the $D_{\widetilde{\mathbf i}}(h')$'s
overlap with which $D_{\mathbf i}(h)$'s. Figure~\ref{f:oldfig}
adapted from \cite{Rie:BruCel}
%\ref{f:oldfig}
shows which $D_{\mathbf i}(h)$'s and $D_{\widetilde{\mathbf
i}}(h')$'s lie in the same connected component, and these
components are numbered for later use. The sequences of signs in
the figure stand for the values of $h:\{1,\dotsc,6\}\to\{\pm
1\}$, and the columns indicate whether $D_{\mathbf i}(h)$ or
$D_{\widetilde{\mathbf i}}(h)$ is meant. So for example the
sequence of signs in the third row first column stands for
$D_{\mathbf i}((1,1,-1,1,-1,1))$, and the figure says that this
cell lies in connected component 5.

%\begin{figure}\label{f:oldfig}
%\epsfbox{G2pic1ij.eps} \epsfbox{G2pic2ij.eps}
%\end{figure}

\section{Berenstein--Zelevinsky's generalized chamber Ansatz for $G_2$}

%our point of view with maps epsilon and alpha, definitions, etc,
%strategy of proof
We wish to determine which connected component of $\mathcal{B}^*$
each $D_{\sigma}(h)$ lies in. We shall do this by proving that
each $D_{\sigma}(h)$ actually intersects one of the cells
$D_{\widetilde{\mathbf i}}(h')$ (for some choice $h'$ of signs)
listed in the previous section. It then follows that
$D_{\sigma}(h)$ lies in the same connected component as this cell.

In order to show that an element $y\in D_{\sigma}^-$ also lies in
one of the cells \linebreak $D_{\widetilde{\mathbf i}}^-(h')$, we
need to show that $y$ can be expressed in the form \linebreak
$y=y_{\widetilde{\mathbf i}}(a_1,\dotsc, a_6) =
y_2(a_1)y_1(a_2)y_2(a_3)y_1(a_4)y_2(a_5)y_1(a_6)$ for some
non--zero real numbers $a_1,a_2,a_3,a_4,a_5,a_6\in \mathbb{R}^*$
with signs $h'(1)$, $h'(2)$, $h'(3)$, $h'(4)$, $h'(5)$, $h'(6)$.
We can then deduce which connected component $D_{\sigma}(h)$ lies
in.

There already exists in the literature a beautiful method for determining
$a_1$, $a_2$, $a_3$, $a_4$, $a_5$, $a_6$,
should they exist --- this is the Chamber Ansatz
of Berenstein, Fomin and Zelevinsky; we therefore employ this method.

%We first of all would like to advance the point of view that Theorem 1.4
%in~\cite{bz4} is actually an explicit calculation of the map
%$\eta:U^+\cap B^-wB^-\rightarrow U^+\cap B^-wB^-$ --- starting
%with $z\in U^+\cap B^+wB^+$, Theorem 1.4 actually describes the
%factorisation of $x=\eta_w(z)\in U^+\cap B^-wB^-$. We first of all give a
%more symmetric description of $\eta_w$, by considering its composition with
%the transpose map. Recall that this is the map corresponding to the
%involutive antiautomorphism of the Lie algebra $\mathbf{g}$ of $G$ given by
%
%$$e_i^T=f_i,\ f_i^T=e_i,\ h_i^T=h_i\ (i=1,\ldots,r);$$
%
%we will use the same notation for the corresponding automorphism of $G$.

There is a natural map $\varepsilon:U^+\cap B^-\dot w_0 B^-
\rightarrow U^-\cap B^+\dot w_0 B^+$, given explicitly by
$\varepsilon= i_-^{-1}\circ i_+$ (see~\ref{e:isos}); similarly we
define $\alpha:U^-\cap B^+wB^+ \rightarrow U^+\cap B^-w^{-1}B^-$
to be $\alpha=i_+^{-1}\circ i_-$, the inverse of $\varepsilon$.

\begin{equation*}
U^+\cap B^-\dot w_0 B^- \underset{\varepsilon} \longrightarrow
U^-\cap B^+\dot w_0 B^+
\end{equation*}
\begin{equation*}
U^+\cap B^-\dot w_0 B^- \underset{\alpha} \longleftarrow
U^-\cap B^+\dot w_0 B^+
\end{equation*}

%If $x\in U^+$, then there is a unique
%element $y\in U^-$ such that $x\cdot B^-=y\cdot U^+$, which we call
%$\epsilon_w(x)$ []reference?Is this correct?. It can be checked that
%$\varepsilon_w=\eta_w^T$ []include?.
%If $y\in U^-$, then there is a unique
%element $x\in U^-$ such that $y\cdot B^+=x\cdot U^-$, which we call
%$\alpha_w(x)$.

Let $y\in U^-\cap B^+\dot w_0 B^+$. Then Theorem 1.4 in~\cite{bz4}
can be used to calculate $\alpha(y)\in U^+\cap B^-\dot w_0B^-$ as
an element of the form
\begin{equation*}
x=x_{\widetilde{\mathbf i}}(a,b,c,d,e,f):=x_2(a) x_1(b) x_2(c)
x_1(d) x_2(e) x_1(f).
\end{equation*}
A second application of this Theorem can be used to calculate
$y=\varepsilon(x)\in U^-\cap B^+\dot w_0B^+$ as an element of the
form $y_{\widetilde{\mathbf i}}(a',b',c',d',e',f')$. If the minors
in the Chamber Ansatz do not vanish, then we have found
$a',b',c',d',e',f'\in \mathbb{R}^*$ such that
$y=y_{\widetilde{\mathbf i}}(a',b',c',d',e',f')$ as required.

We shall therefore need a description of the Chamber Ansatz in
type $G_2$. If $w\in W$ and $i=1,2$, then Berenstein and
Zelevinsky define $w\omega_i$ to be a {\em chamber weight of
level $i$}, and define a corresponding `generalised minor'
$\Delta^{w\omega_i}$, which is a function on $G$. These minors
reduce to usual minors of a matrix in type $A$. The generalised
minors can be described (using section 6 of~\cite{bz4}) in the
following way. Let $g\in G$, $w\in W$, and $\omega$ a dominant
weight. Then there is a module $V_{\omega}$ for $G$ of highest
weight $\omega$. Fix a vector $v_{\omega}\in V_{\omega}$ of
highest weight $\omega$, fix a reduced expression
$w=s_{j_1}s_{j_2}\cdots s_{j_l}$ for $w$, and for $k=1,\ldots ,l$
set $b_k=\left <\alpha^{\vee}_{j_k},s_{j_{k-1}}\cdots
s_{j_1}\omega\right  >$. Then
$$v_{w\omega}:=f_{j_l}^{(b_l)}f_{j_{l-1}}^{(b_{(l-1)})}\cdots
f_{j_1}^{(b_1)}v_{\omega}$$ is an extremal weight vector of
weight $w\omega$. Then $\Delta^{w\omega}(g)$ is defined to be the
coefficient of $v_{\omega}$ in $g\cdot v_{w\omega}$. We have:

\begin{thm} (Berenstein and Zelevinsky) \label{ansatz1}
Let $x\in U^+\cap B^-\dot w_0 B^-$, and suppose that
$\mathbf{j}=(j_1,\dotsc , j_N)$ is a reduced expression for $w_0$.
Then
\begin{equation*}
\varepsilon(x)=y=y_{j_1}(a_1)y_{j_2}(a_2)\cdots y_{j_N}(a_N),
\end{equation*}
where $a_1,a_2,\ldots ,a_N$ are given by:
$$a_{N+1-k}=\frac{1}{\Delta^{w_k\omega_{j_k}}(x)\Delta^{w_{k+1}\omega_{j_k}}(x)}
\prod_{j\not=j_k}\Delta^{w_k\omega_j}(x)^{-A_{j,j_k}},$$

where $w_k=s_{j_N}s_{j_{N-1}}\cdots s_{j_k}$.
\end{thm}

If we are in the situation of the Theorem, we write
$\varepsilon^{\mathbf{j}}(x)=(a_1,a_2,\ldots ,a_N)$.

Again using~\cite[Theorem 1.4]{bz4}, we can describe $\alpha$ using the
Chamber Ansatz. Define
$\Delta_-^{w\omega}(g)$ to be the coefficient of $v_{-\omega}$ in
$g\cdot v_{w\omega}$, where $v_{w\omega}$ is the extremal weight vector
of weight $w\omega$ as defined above.

\begin{thm} (Berenstein and Zelevinsky) \label{ansatz2}
Let $y\in U^-\cap B^+\dot w_0 B^+$, and suppose that
$\mathbf{j}=(j_1,\dotsc , j_N)$ is a reduced expression for $w_0$.
Then
\begin{equation*}
\alpha(y)=x=x_{j_1}(a_1)x_{j_2}(a_2)\cdots x_{j_N}(a_N),
\end{equation*}
where $a_1,a_2,\ldots ,a_N$ are given by:

$$a_{N+1-k}=\frac{1}{\Delta_-^{-w_k\omega_{j_k}}(y)\Delta_-^{-w_{k+1}\omega_{j_k}}(y)}
\prod_{j\not=j_k}\Delta_-^{-w_k\omega_j}(y)^{-A_{j,j_k}},$$

where $w_k=s_{j_N}s_{j_{N-1}}\cdots s_{j_k}$.
\end{thm}

If we are in the situation of this Theorem, we write
$\alpha^{\mathbf{j}}(x)=(a_1,a_2,\ldots ,a_N)$.

\begin{figure} \label{chamberdiagram1}
\epsfbox{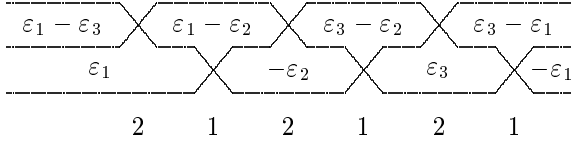} \caption{The Chamber Ansatz for $\varepsilon^{212121}$.}
\end{figure}

We can use a wiring diagram to describe these maps in type $G_2$,
much as in type $A$. Recall that our two reduced expressions of
$w_0$ are $\mathbf i=(1,2,1,2,1,2)$ and $\widetilde{\mathbf
i}=(2,1,2,1,2,1)$. We begin by calculating need
$\varepsilon^{212121}=\varepsilon^{\widetilde{\mathbf i}}$; this
can be done using the Chamber Ansatz in
Figure~\ref{chamberdiagram1}.

Suppose that $a_j$ is one of the parameters in
Theorem~\ref{ansatz1}, and suppose that adjacent to the $j$th
crossing from the left in Figure~\ref{chamberdiagram1}, chamber
$A$ is above, $D$ below, and $B$ and $C$ on the same horizontal
level. Let $\Delta(A)$, $\Delta(B)$, $\Delta(C)$, $\Delta(D)$ be
the generalized minors corresponding to these four chambers
(taking the value $1$ on unbounded chambers). Then
Theorem~\ref{ansatz1} states that:

$$a_j=\frac{\Delta(A)^3\Delta(D)^3}{\Delta(B)\Delta(C)}.$$

\begin{figure} \label{chamberdiagram2}
\epsfbox{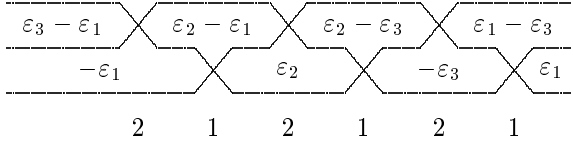} \caption{The Chamber Ansatz for $\alpha^{212121}$.}
\end{figure}

We also need to calculate $\alpha^{\widetilde{\mathbf
i}}=\alpha^{212121}$, which can be done using the Chamber Ansatz
in Figure~\ref{chamberdiagram2}. The value of the component at a
crossing is calculated as for $\varepsilon$, but using this
second picture.

\section{The map $i_+$}

\begin{center}
\begin{figure}
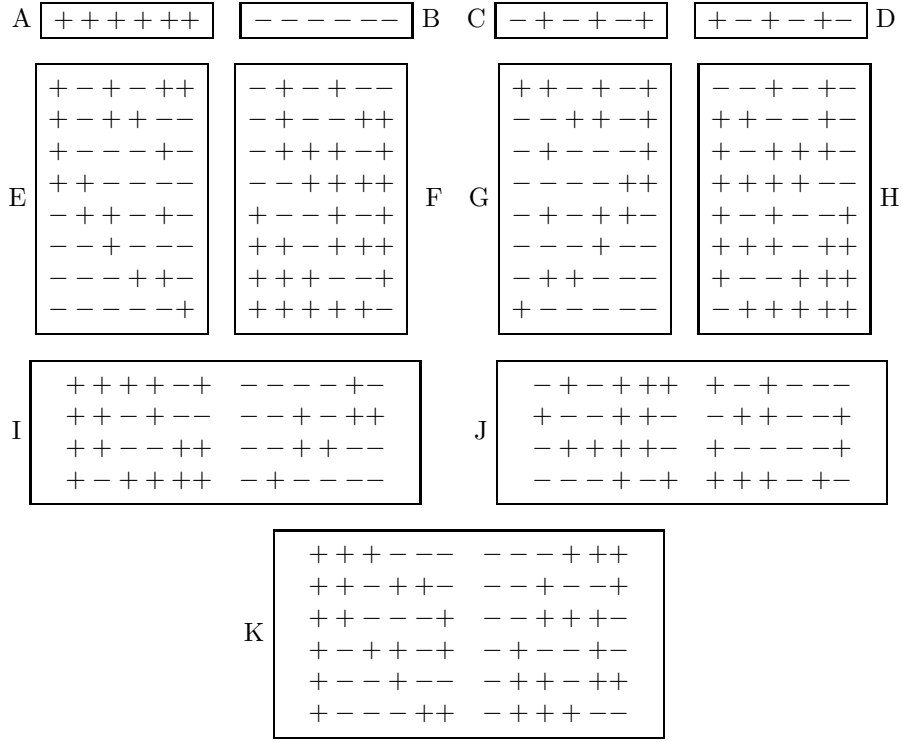

\begin{tabular}{cccc}
A\ \framebox[2.3cm]{$++++++$}&\framebox[2.3cm]{$------$}\ B & C\
\framebox[2.3cm]{$-+-+-+$}& \framebox[2.3cm]{$+-+-+-$}\ D
\end{tabular}

\vspace*{0.3cm}

\begin{tabular}{cccc}
E\
 \framebox[2.3cm]{
$
\begin{array}{cc}
+-+-++  \\
+-++-- \\
+---+- \\
++---- \\
-++-+- \\
--+--- \\
---++- \\
-----+
\end{array}
$ } &
 \framebox[2.3cm]{
$
\begin{array}{cc}
-+-+-- \\
-+--++\\
-+++-+\\
--++++\\
+--+-+\\
++-+++\\
+++--+\\
+++++-
\end{array}$ }
\ F & G\
 \framebox[2.3cm]{
$
\begin{array}{cc}
++-+-+\\
--++-+\\
-+---+\\
----++\\
-+-++-\\
---+--\\
-++---\\
+-----
\end{array}$ } &
 \framebox[2.3cm]{
$
\begin{array}{cc}
--+-+-\\
++--+-\\
+-+++-\\
++++--\\
+-+--+\\
+++-++\\
+--+++\\
-+++++\\
\end{array}$ }\ H
\end{tabular}

\vspace*{0.3cm}

\begin{tabular}{cccc}
\ \ I\  \framebox[5.2cm]{ $
\begin{array}{cc}
 ++++-+ & ----+-\\
 ++-+-- & --+-++\\
 ++--++ & --++--\\
 +-++++ & -+----
\end{array}$ } &
\quad J\ \framebox[5.2cm]{ $
\begin{array}{cc}
 -+-+++ & +-+---\\
 +--++- & -++--+\\
 -++++- & +----+\\
 ---+-+ & +++-+-
\end{array}$ } &
\end{tabular}

\vspace*{0.3cm}

\begin{tabular}{ccc}
&K\ \framebox[5.2cm]{ $
\begin{array}{cc}
 +++--- & ---+++ \\
 ++-++- & --+--+ \\
 ++---+ & --+++- \\
 +-++-+ & -+--+- \\
 +--+-- & -++-++ \\
 +---++ & -+++--
\end{array}$ }&
\end{tabular}
\caption{\label{f:varfig} The connected components of $U^+\cap
B^-\dot w_0 B^-$ in terms of the cells $D_{\widetilde{\mathbf
i}}^+(h)$.}
\end{figure}
\end{center}

 Note that our approach to parameterizing $\mathcal B^*$ in
Section \ref{s:param} was inherently asymmetric, suited to the
map $i_-:U^-\cap B^+\dot w_0 B^+ \to \mathcal B^*$. That is, the
cells $D_{\mathbf j}(h)$ came naturally from cells $D^-_{\mathbf
j}(h)$ defined in $U^-$. We now need to consider another set of
cells in $\mathcal B^*$ coming from the isomorphism
 $
 i^+:U^+\cap B^- \dot w_0 B^-\to \mathcal B^*.
 $
Let us define
\begin{eqnarray*}
x_{\widetilde{\mathbf i}}(a_1,\dotsc,
a_6)&:=& x_2(a_1) x_1(a_2) x_2(a_3)
x_1(a_4) x_2(a_5) x_1(a_6),\\
D^+_{\widetilde{\mathbf i}}(h)&:=& \{x_{\widetilde{\mathbf i}}(a_1,\dotsc,
a_6)\ |\  a_i h(i)\in\R_{>0}\}.
\end{eqnarray*}
%And set $\mathcal D^+_{\widetilde{\mathbf
%i}}(h):=i_+(D^+_{\widetilde{\mathbf i}}(h))$, the corresponding
%cell in $\mathcal B^*$.
for any map $h:\{1,2,3,4,5,6\}\to\{\pm 1\}$. It is clear that two
cells $D^+_{\widetilde{\mathbf i}}(h)$ and
$D^+_{\widetilde{\mathbf i}}(h')$ lie in the same connected
component of $U^+\cap B^+\dot w_0 B^+$ precisely if
$D^-_{\widetilde{\mathbf i}}(h)$ and $D^-_{\widetilde{\mathbf
i}}(h')$ lie in the same component, by symmetry
between $U^+$ and $U^-$. Reading off
from Figure \ref{f:oldfig} when this happens, we obtain the new
Figure \ref{f:varfig}, which groups together all sequences of
signs (determining maps $h:\{1,2,3,4,5,6\}\to\{\pm 1\}$) such
that the corresponding $D^+_{\widetilde{\mathbf i}}(h)$'s lie in
the same connected component of $U^+\cap B^-\dot w_0 B^-$. We
have labeled these connected components by letters $A$--$K$.

Our aim is now to identify which connected components of $U^+\cap
B^-\dot w_0 B^-$ labeled by $A$--$K$ correspond to which
components labeled $1$--$11$ of $\mathcal B^*$ under the map
$i_+$. To do this we simply need to evaluate the map
$\varepsilon^{212121}$ on $11$ test points, one from each
component of $U^+\cap B^-\dot w_0 B^-$.

It suffices to calculate $\varepsilon^{212121}$ on elements $x$ of
the form $x=x_{\widetilde{\mathbf i}}(a,b,c,d,e,f)$.
%It is easy to see that $\varepsilon$ maps a cell $D_{\mathbf
%i}^+(h)$ of $U^+\cap B^-\dot w_0 B^-$ to inside a connected
%component of $U^-\cap B^+\dot w_0 B^+$. We know what these are
%(see Section \ref{s:param}). We can therefore determine which
%component this is by substituting in appropriate values for
%$a,b,c,d,e,f$ and then applying $\varepsilon^{212121}$; the signs
%of the new parameters obtained determine which component
%$\varepsilon(x)$ lies in.
%First of all, we calculate:
%$x_{212121}(a,b,c,d,e,f)=$
%$$\tiny
%\left [\begin {array}{ccccccc} 1&e+c+a&\left (e+c+a\right )f+\left (c+
%a\right )d+ab&2\,\left (\left (c+a\right )d+ab\right )e+2\,abc&\left (
%\left (c+a\right )d+ab\right ){e}^{2}+2\,abce+ab{c}^{2}&\left (\left (
%\left (c+a\right )d+ab\right ){e}^{2}+2\,abce+ab{c}^{2}\right )f+ab{c}
%^{2}d&ab{c}^{2}de\\\noalign{\medskip}0&1&f+d+b&2\,\left (d+b\right )e+
%2\,bc&\left (d+b\right ){e}^{2}+2\,bce+b{c}^{2}&\left (\left (d+b
%\right ){e}^{2}+2\,bce+b{c}^{2}\right )f+b{c}^{2}d&b{c}^{2}de
%\\\noalign{\medskip}0&0&1&2\,e+2\,c+2\,a&{e}^{2}+\left (2\,c+2\,a
%\right )e+{c}^{2}+2\,ac+{a}^{2}&\left ({e}^{2}+\left (2\,c+2\,a\right
%)e+{c}^{2}+2\,ac+{a}^{2}\right )f+\left ({c}^{2}+2\,ac+{a}^{2}\right )
%d+{a}^{2}b&\left (\left ({c}^{2}+2\,ac+{a}^{2}\right )d+{a}^{2}b
%\right )e+{a}^{2}bc\\\noalign{\medskip}0&0&0&1&e+c+a&\left (e+c+a
%\right )f+\left (c+a\right )d+ab&\left (\left (c+a\right )d+ab\right )
%e+abc\\\noalign{\medskip}0&0&0&0&1&f+d+b&\left (d+b\right )e+bc
%\\\noalign{\medskip}0&0&0&0&0&1&e+c+a\\\noalign{\medskip}0&0&0&0&0&0&1
%\end {array}\right ].
%$$
%%Maple did this!
From the Chamber Ansatz we have
\begin{equation*} \varepsilon(x_{\widetilde{\mathbf
i}}(a,b,c,d,e,f))=y_{\widetilde{\mathbf i}}(a',b',c',d',e',f'),
\end{equation*}
where
\begin{eqnarray*}
a' & = & \frac{1}{\Delta^{\varepsilon_1-\varepsilon_3}(x)\Delta^{\varepsilon_1-\varepsilon_2}(x)}, \\
b' & = & \frac{\Delta^{\varepsilon_1-\varepsilon_2}(x)}{\Delta^{\varepsilon_1}(x)\Delta^{-\varepsilon_2}(x)}, \\
c' & = & \frac{\Delta^{-\varepsilon_2}(x)^3}{\Delta^{\varepsilon_1-\varepsilon_2}(x)\Delta^{\varepsilon_3-\varepsilon_2(x)}}, \\
d' & = & \frac{\Delta^{\varepsilon_3-\varepsilon_2}(x)}{\Delta^{-\varepsilon_2}(x)\Delta^{\varepsilon_3}(x)}, \\
e' & = & \frac{\Delta^{\varepsilon_3}(x)^3}{\Delta^{\varepsilon_3-\varepsilon_2}(x)\Delta^{\varepsilon_3-\varepsilon_1}(x)}, \\
f' & = &
\frac{\Delta^{\varepsilon_3-\varepsilon_1}(x)}{\Delta^{\varepsilon_3}(x)\Delta^{-\varepsilon_1}(x)}.
\end{eqnarray*}

%In type $G_2$, there are two fundamental representations:
%$V_{\omega_1}$, the $7$--dimensional representation, and
%$V_{\omega_2}$, the $14$--dimensional adjoint representation.
By considering the usual action of the Chevalley generators of the
Lie algebra in the two fundamental representations $V_{\omega_1}$
and $V_{\omega_2}$ of $G_2$, the action of
$x_1(t)$,$x_2(t)$,\linebreak
$y_1(t)$,$y_2(t)$,$\dot s_1$ and $\dot s_2$ can be
computed explicitly. We list below the relevant matrix
coefficients for the element $x=x_{\widetilde{\mathbf
i}}(a,b,c,d,e,f)$.

%We obtain:

\begin{eqnarray*}
  \Delta^{\varepsilon_1}(x) & = & 1, \\
  \Delta^{-\varepsilon_3}(x) & = & f+d+b, \\
  \Delta^{-\varepsilon_2}(x) & = & ed+eb+bc, \\
  \Delta^{\varepsilon_2}(x) & = & {f}^{2}ed+{f}^{2}eb+{f}^{2}bc+2\,bcdf+bc{d}^{2}, \\
  \Delta^{\varepsilon_3}(x) & = & bc{d}^{2}e, \\
  \Delta^{-\varepsilon_1}(x) & = & bc{d}^{2}ef, \\
 \Delta^{\varepsilon_1-\varepsilon_3}(x) & = & 1, \\
%\end{eqnarray*}
%\begin{eqnarray*}
 \Delta^{\varepsilon_1-\varepsilon_2}(x) & = & e+c+a, \\
\Delta^{\varepsilon_2-\varepsilon_3}(x)
 & = & {f}^{3}e+{f}^{3}c+{f}^{3}a+3\,{f}^{2}cd+3\,{f}^{2}ad+3\,{f}^{2}ab+3\,f
{d}^{2}c+3\,f{d}^{2}a+\, \\
 \lefteqn{{}+6fabd+3\,fa{b}^{2}+{d}^{3}c+{d}^{3}a+3\,ab{d}^{2}+3\,a{b}^{2}d+a{b}^{3},} \\
 \Delta^{\varepsilon_3-\varepsilon_2}(x) & = &
{e}^{2}{d}^{3}c+{e}^{2}{d}^{3}a+3\,{e}^{2}ab{d}^{2}+3\,{e}^{2}a{b}^{2}
d+{e}^{2}a{b}^{3}+3\,ea{b}^{2}cd+2\,ea{b}^{3}c \\
\lefteqn{{}+a{b}^{3}{c}^{2},} \\
 \Delta^{\varepsilon_2-\varepsilon_1}(x) & = &
{f}^{3}{e}^{2}{d}^{3}c+{f}^{3}{e}^{2}{d}^{3}a+3\,{f}^{3}{e}^{2}ab{d}^{
2}+3\,{f}^{3}{e}^{2}a{b}^{2}d+{f}^{3}{e}^{2}a{b}^{3}+ \\
\lefteqn{{}+3\,{f}^{3}ea{b}^{
2}cd+2\,{f}^{3}ea{b}^{3}c+{f}^{3}a{b}^{3}{c}^{2}+3\,{f}^{2}ea{b}^{2}c{
d}^{2}+3\,{f}^{2}ea{b}^{3}cd+} \\
\lefteqn{{}+3\,{f}^{2}a{b}^{3}{c}^{2}d+3\,a{b}^{3}{c}
^{2}{d}^{2}f+a{b}^{3}{c}^{2}{d}^{3},} \\
 \Delta^{\varepsilon_3-\varepsilon_1}(x) & = &
a{b}^{3}{c}^{2}{d}^{3}e.
\end{eqnarray*}

Applying the Chamber Ansatz, we obtain the components $a',b',c',d',e',f'$:

\begin{eqnarray*}
   a' & = &
\frac{1}{e+c+a}, \\
   b' & = &
{\frac {e+c+a}{ed+eb+bc}}, \\
   c' & = &
{\frac {\left (ed+eb+bc\right )^{3}}{u\left (e+c+a\right )}}, \\
   d' & = &
{\frac {u}{bc{d}^{2}e\left (ed+eb+bc\right )}}, \\
   e' & = &
{\frac {{e}^{2}{d}^{3}c}{au}}, \\
   f' & = &
{\frac {ab}{def}},
\end{eqnarray*}

where
 $$
 u={\left ({e}^{2}{d}^{3}c+{e}^{2}{d}
^{3}a+3\,{e}^{2}ab{d}^{2}+3\,{e}^{2}a{b}^{2}d+{e}^{2}a{b}^{3}+3\,ea{b}
^{2}cd+2\,ea{b}^{3}c+a{b}^{3}{c}^{2}\right )}.
 $$
The next step is to substitute values for $a,b,c,d,e,f$ with
specific combinations of signs, one choice from each component of
Figure~\ref{f:varfig}, ensuring that the matrix coefficients in
the Chamber Ansatz do not vanish. It turns out that this can be
achieved by setting $a=\pm 1,b=\pm 2, c=\pm 3, d=\pm 5, e=\pm 7$
and $f=\pm 11$. For any choice of signs of $a,b,c,d,e,f$, we
obtain the signs of $a',b',c',d',e',f'$, and thus the connected
component using Figure~\ref{f:oldfig}. The resulting bijection
between connected components of $U^+\cap B^-\dot w_0 B^-$ and
$\mathcal B^*$ is shown below.
\begin{equation}\label{e:bij}
\begin{array}{c}
\begin{array}{cccc}
A\longleftrightarrow 1& B\longleftrightarrow
2&C\longleftrightarrow 3&D\longleftrightarrow 4\\
E\longleftrightarrow 5&F\longleftrightarrow
6&G\longleftrightarrow 7&H\longleftrightarrow 8\\
\end{array}\\
\begin{array}{cc}
I\longleftrightarrow 10& J\longleftrightarrow 9
\end{array}\\
K\longleftrightarrow 11
\end{array}
\end{equation}

%\begin{center} Table I \\ Bijection of connected components under $i_+$ or equivalently $\epsilon$ \\ \ \\
%\begin{table}{c|c}
%Connected component of $U^+\cap B^-\dot w_0 B^-$ & Connected component of $\mathcal B^*$ \\
%\hline
%\endfirsthead
%Signs of $a,b,c,d,e,f$ &  Connected component \\
%\hline
%\endhead
%\endfoot
%\endlastfoot
%  A &  1 \\
%  B &  2 \\
%  C &  3 \\
%  D &  4 \\
%  E &  5\\
%  F &  6\\
%  G &  7\\
%  H &  8 \\
%  I &  10\\
%  J &  9\\
%  K &  11\\
%\end{table}
%\end{center}

%Which
%connected component can be read off Figure~\ref{f:oldfig}; they
%are numbered as shown there, starting at the top and moving from
%left to right, down the list, giving components 1--11.

%\vspace*{0.3cm}

\section{Calculation of the Euler Characteristics}

We need to determine which connected components of $U^-\cap
B^+\dot w_0B^+$ the $D_{\sigma}^-(h)$ belong to. If one point in
such a $D_{\sigma}^-(h)$ lies in a particular connected
component, then the whole of it does. So a similar approach to
that in the previous section will work. We start with a general
point in $D_{\sigma}^-(h)\subseteq U^-\cap B^+\dot w_0  B^+$,
apply $\alpha$ to it, to get a point in $U^+\cap B^-\dot w_0
B^-$; we ensure this is of the form
$x_{212121}(a',b',c',d',e',f')$ for some $a',b',c',d',e',f'$;
then it lies in one of the components $A$--$K$ of
Figure~\ref{f:varfig}; the bijection \eqref{e:bij} above will then
give us the connected component of $D_{\sigma}^-(h)$. Again, not
every point in $D_{\sigma}^-$ need be in the domain of
$\alpha^{212121}$ but we can always find one which is. We use the
Chamber Ansatz for $\alpha^{212121}$ as described above. Our
computation is as follows:

\begin{center}
The following are the components of $\alpha^{212121}$ applied to
$y_{x21x12}(t_1,t_2,m_1,m_2)$ (provided this is well-defined).

\begin{eqnarray*}
a & = &
-{{\it m_2}}^{-1}, \\
b & = &
{\frac {{\it m_2}}{{\it t_1}\,{\it m_2}+{\it m_1}}}, \\
c & = &
{\frac {\left ({\it t_1}\,{\it m_2}+{\it m_1}\right )^{3}}{{\it m_2}\,
\left (-{\it m_2}\,{\it t_2}+{{\it m_1}}^{3}\right )}}, \\
d & = &
{\frac {-{\it m_2}\,{\it t_2}+{{\it m_1}}^{3}}{\left ({\it t_1}\,{\it m_2}+
{\it m_1}\right )\left ({\it t_1}\,{{\it m_1}}^{2}+{\it t_2}\right )}}, \\
e & = &
-{\frac {\left ({\it t_1}\,{{\it m_1}}^{2}+{\it t_2}\right )^{3}}{\left (
-{\it m_2}\,{\it t_2}+{{\it m_1}}^{3}\right ){{\it t_2}}^{2}}}, \\
f & = &
{\frac {{\it t_2}}{\left ({\it t_1}\,{{\it m_1}}^{2}+{\it t_2}\right ){
\it t_1}}}.
\end{eqnarray*}
\end{center}

\begin{center}
The following are the components of $\alpha^{212121}$ applied to
$y_{12x21x}(t_1,m_1,m_2, t_2)$.

\begin{eqnarray*}
a & = &
{{\it t_2}}^{-1}, \\
b & = &
-{{\it m_2}}^{-1}, \\
c & = &
{\frac {{{\it m_2}}^{3}}{3\,{\it t_1}\,{\it m_2}+{\it m_1}}}, \\
d & = &
{\frac {3\,{\it t_1}\,{\it m_2}+{\it m_1}}{{\it m_2}\,\left (2\,{\it t_1}\,
{\it m_2}+{\it m_1}\right )}}, \\
e & = &
{\frac {\left (2\,{\it t_1}\,{\it m_2}+{\it m_1}\right )^{3}}{\left (3\,{
\it t_1}\,{\it m_2}+{\it m_1}\right ){{\it t_1}}^{3}}}, \\
f & = &
-{\frac {{\it t_1}}{2\,{\it t_1}\,{\it m_2}+{\it m_1}}}.
\end{eqnarray*}
\end{center}

\begin{center}
The following are the components of $\alpha^{212121}$ applied to
$y_{1x12x2}(t_1,m_1,t_2,m_2)$.

\begin{eqnarray*}
a & = &
-{{\it m_2}}^{-1}, \\
b & = &
-{\frac {{\it m_2}}{{\it m_1}\,{\it m_2}+{\it t_2}}}, \\
c & = &
-{\frac {\left ({\it m_1}\,{\it m_2}+{\it t_2}\right )^{3}}{{\it m_2}\,
\left ({\it t_1}\,{{\it m_2}}^{2}-{{\it t_2}}^{3}\right )}}, \\
d & = &
{\frac {{\it t_1}\,{{\it m_2}}^{2}-{{\it t_2}}^{3}}{\left ({\it m_1}\,{
\it m_2}+{\it t_2}\right )\left ({\it t_1}\,{\it m_2}+{\it m_1}\,{{\it t_2}}
^{2}\right )}}, \\
e & = &
{\frac {\left ({\it t_1}\,{\it m_2}+{\it m_1}\,{{\it t_2}}^{2}\right )^{3}
}{\left ({\it t_1}\,{{\it m_2}}^{2}-{{\it t_2}}^{3}\right ){\it t_1}\,{{
\it t_2}}^{3}}}, \\
f & = &
{\frac {{{\it t_2}}^{2}}{{\it t_1}\,{\it m_2}+{\it m_1}\,{{\it t_2}}^{2}}}.
\end{eqnarray*}
\end{center}

\begin{center}
The following are the components of $\alpha^{212121}$ applied to
$y_{xx1x1x}(t_1,t_2,t_3,m_1, t_4)$.

\begin{eqnarray*}
a & = &
\left ({\it t_2}+{\it t_4}\right )^{-1}, \\
b & = &
{\frac {{\it t_2}+{\it t_4}}{{\it t_1}\,{\it t_2}-{\it m_1}\,{\it t_4}+{\it
t_4}\,{\it t_1}}}, \\
c & = &
-{\frac {\left ({\it t_1}\,{\it t_2}-{\it m_1}\,{\it t_4}+{\it t_4}\,{\it
t_1}\right )^{3}}{\left ({\it t_2}+{\it t_4}\right ){\it t_4}\,\left (-{
\it t_2}\,{\it t_3}+{\it t_4}\,{{\it m_1}}^{3}{\it t_2}-{\it t_4}\,{\it t_3}
\right )}}, \\
d & = &
-{\frac {-{\it t_2}\,{\it t_3}+{\it t_4}\,{{\it m_1}}^{3}{\it t_2}-{\it t_4}
\,{\it t_3}}{\left ({\it t_1}\,{\it t_2}\,{{\it m_1}}^{2}-{\it t_3}\right )
\left ({\it t_1}\,{\it t_2}-{\it m_1}\,{\it t_4}+{\it t_4}\,{\it t_1}\right
)}}, \\
e & = &
{\frac {\left ({\it t_1}\,{\it t_2}\,{{\it m_1}}^{2}-{\it t_3}\right )^{3}
{\it t_4}}{\left (-{\it t_2}\,{\it t_3}+{\it t_4}\,{{\it m_1}}^{3}{\it t_2}-
{\it t_4}\,{\it t_3}\right ){\it t_2}\,{{\it t_3}}^{2}}}, \\
f & = &
-{\frac {{\it t_3}}{\left ({\it t_1}\,{\it t_2}\,{{\it m_1}}^{2}-{\it t_3}
\right ){\it t_1}}}.
\end{eqnarray*}
\end{center}

\begin{center}
The following are the components of $\alpha^{212121}$ applied to
$y_{1x1xxx}(t_1,m_1,t_2,t_3,t_4)$.

\begin{eqnarray*}
a & = &
\left ({\it t_2}+{\it t_4}\right )^{-1}, \\
b & = &
-{\frac {{\it t_2}+{\it t_4}}{{\it m_1}\,{\it t_2}+{\it m_1}\,{\it t_4}-{
\it t_4}\,{\it t_3}}}, \\
c & = &
-{\frac {\left ({\it m_1}\,{\it t_2}+{\it m_1}\,{\it t_4}-{\it t_4}\,{\it
t_3}\right )^{3}}{\left ({\it t_2}+{\it t_4}\right )\left ({\it t_1}\,{{
\it t_2}}^{2}+2\,{\it t_1}\,{\it t_2}\,{\it t_4}+{{\it t_4}}^{2}{\it t_1}+{
\it t_2}\,{{\it t_3}}^{3}{{\it t_4}}^{2}\right )}}, \\
d & = &
{\frac {{\it t_1}\,{{\it t_2}}^{2}+2\,{\it t_1}\,{\it t_2}\,{\it t_4}+{{
\it t_4}}^{2}{\it t_1}+{\it t_2}\,{{\it t_3}}^{3}{{\it t_4}}^{2}}{\left ({
\it m_1}\,{\it t_2}+{\it m_1}\,{\it t_4}-{\it t_4}\,{\it t_3}\right )\left (
{\it t_1}\,{\it t_2}+{\it t_4}\,{\it t_1}+{\it t_2}\,{{\it t_3}}^{2}{\it t_4}
\,{\it m_1}\right )}}, \\
e & = &
-{\frac {\left ({\it t_1}\,{\it t_2}+{\it t_4}\,{\it t_1}+{\it t_2}\,{{\it
t_3}}^{2}{\it t_4}\,{\it m_1}\right )^{3}}{\left ({\it t_1}\,{{\it t_2}}^{2
}+2\,{\it t_1}\,{\it t_2}\,{\it t_4}+{{\it t_4}}^{2}{\it t_1}+{\it t_2}\,{{
\it t_3}}^{3}{{\it t_4}}^{2}\right ){\it t_1}\,{{\it t_2}}^{2}{{\it t_3}}^{
3}{\it t_4}}}, \\
f & = &
{\frac {{\it t_2}\,{{\it t_3}}^{2}{\it t_4}}{{\it t_1}\,{\it t_2}+{\it t_4}
\,{\it t_1}+{\it t_2}\,{{\it t_3}}^{2}{\it t_4}\,{\it m_1}}}.
\end{eqnarray*}
\end{center}

\begin{center}
The following are the components of $\alpha^{212121}$ applied to
$y_{xxx2x2}(t_1,t_2,t_3,t_4,m_1)$.

\begin{eqnarray*}
a & = &
-\left ({\it m_1}-{\it t_2}\right )^{-1}, \\
b & = &
{\frac {{\it m_1}-{\it t_2}}{{\it t_1}\,{\it m_1}+{\it m_1}\,{\it t_3}-{\it
t_1}\,{\it t_2}-{\it t_4}}}, \\
c & = &
{\frac {\left ({\it t_1}\,{\it m_1}+{\it m_1}\,{\it t_3}-{\it t_1}\,{\it t_2
}-{\it t_4}\right )^{3}}{\left ({\it m_1}-{\it t_2}\right )\left ({\it t_2
}\,{{\it t_3}}^{3}{{\it m_1}}^{2}-3\,{\it t_2}\,{{\it t_3}}^{2}{\it t_4}\,{
\it m_1}+3\,{\it t_2}\,{\it t_3}\,{{\it t_4}}^{2}-{{\it t_4}}^{3}\right )}}, \\
d & = &
{\frac {{\it t_2}\,{{\it t_3}}^{3}{{\it m_1}}^{2}-3\,{\it t_2}\,{{\it t_3}}
^{2}{\it t_4}\,{\it m_1}+3\,{\it t_2}\,{\it t_3}\,{{\it t_4}}^{2}-{{\it t_4}
}^{3}}{\left ({\it t_1}\,{\it m_1}+{\it m_1}\,{\it t_3}-{\it t_1}\,{\it t_2}
-{\it t_4}\right )\left ({\it t_1}\,{\it t_2}\,{{\it t_3}}^{2}{\it m_1}-2\,
{\it t_1}\,{\it t_2}\,{\it t_3}\,{\it t_4}+{{\it t_4}}^{2}{\it t_1}+{\it t_3}
\,{{\it t_4}}^{2}\right )}}, \\
e & = &
-{\frac {\left ({\it t_1}\,{\it t_2}\,{{\it t_3}}^{2}{\it m_1}-2\,{\it t_1}
\,{\it t_2}\,{\it t_3}\,{\it t_4}+{{\it t_4}}^{2}{\it t_1}+{\it t_3}\,{{\it
t_4}}^{2}\right )^{3}}{\left ({\it t_2}\,{{\it t_3}}^{3}{{\it m_1}}^{2}-3
\,{\it t_2}\,{{\it t_3}}^{2}{\it t_4}\,{\it m_1}+3\,{\it t_2}\,{\it t_3}\,{{
\it t_4}}^{2}-{{\it t_4}}^{3}\right ){\it t_2}\,{{\it t_3}}^{3}{{\it t_4}}^
{3}}}, \\
f & = &
{\frac {{\it t_3}\,{{\it t_4}}^{2}}{\left ({\it t_1}\,{\it t_2}\,{{\it t_3}
}^{2}{\it m_1}-2\,{\it t_1}\,{\it t_2}\,{\it t_3}\,{\it t_4}+{{\it t_4}}^{2}
{\it t_1}+{\it t_3}\,{{\it t_4}}^{2}\right ){\it t_1}}}.
\end{eqnarray*}
\end{center}

\begin{center}
The following are the components of $\alpha^{212121}$ applied to
$y_{x2x2xx}(t_1,t_2,m_1,t_3,t_4)$.

\begin{eqnarray*}
a & = &
-\left ({\it m_1}-{\it t_4}\right )^{-1}, \\
b & = &
{\frac {{\it m_1}-{\it t_4}}{{\it t_1}\,{\it m_1}-{\it t_2}-{\it t_4}\,{\it
t_1}-{\it t_4}\,{\it t_3}}}, \\
c & = &
-{\frac {\left ({\it t_1}\,{\it m_1}-{\it t_2}-{\it t_4}\,{\it t_1}-{\it t_4
}\,{\it t_3}\right )^{3}}{\left ({\it m_1}-{\it t_4}\right )\left ({{\it
t_2}}^{3}+3\,{{\it t_2}}^{2}{\it t_3}\,{\it t_4}+3\,{\it t_2}\,{{\it t_3}}^{
2}{{\it t_4}}^{2}+{{\it t_4}}^{2}{\it m_1}\,{{\it t_3}}^{3}\right )}}, \\
d & = &
-{\frac {{{\it t_2}}^{3}+3\,{{\it t_2}}^{2}{\it t_3}\,{\it t_4}+3\,{\it t_2
}\,{{\it t_3}}^{2}{{\it t_4}}^{2}+{{\it t_4}}^{2}{\it m_1}\,{{\it t_3}}^{3}
}{\left ({\it t_1}\,{\it m_1}-{\it t_2}-{\it t_4}\,{\it t_1}-{\it t_4}\,{
\it t_3}\right )\left ({\it t_1}\,{{\it t_2}}^{2}+2\,{\it t_1}\,{\it t_2}\,
{\it t_3}\,{\it t_4}+{\it t_4}\,{{\it t_3}}^{2}{\it t_1}\,{\it m_1}-{\it t_2}
\,{{\it t_3}}^{2}{\it t_4}\right )}}, \\
e & = &
-{\frac {\left ({\it t_1}\,{{\it t_2}}^{2}+2\,{\it t_1}\,{\it t_2}\,{\it
t_3}\,{\it t_4}+{\it t_4}\,{{\it t_3}}^{2}{\it t_1}\,{\it m_1}-{\it t_2}\,{{
\it t_3}}^{2}{\it t_4}\right )^{3}}{\left ({{\it t_2}}^{3}+3\,{{\it t_2}}^
{2}{\it t_3}\,{\it t_4}+3\,{\it t_2}\,{{\it t_3}}^{2}{{\it t_4}}^{2}+{{\it
t_4}}^{2}{\it m_1}\,{{\it t_3}}^{3}\right ){{\it t_2}}^{3}{{\it t_3}}^{3}{
\it t_4}}}, \\
f & = &
-{\frac {{\it t_2}\,{{\it t_3}}^{2}{\it t_4}}{\left ({\it t_1}\,{{\it t_2}}
^{2}+2\,{\it t_1}\,{\it t_2}\,{\it t_3}\,{\it t_4}+{\it t_4}\,{{\it t_3}}^{2
}{\it t_1}\,{\it m_1}-{\it t_2}\,{{\it t_3}}^{2}{\it t_4}\right ){\it t_1}}}.
\end{eqnarray*}
\end{center}

We now substitute in values for the $t_i$ and the $m_i$ in order to obtain
points in the $D_{\sigma}^-(h)$ on which $\alpha^{212121}$ is defined.
In the first three cases we set
$t_1=\pm 1$, $t_2=\pm 2$, $m_1=3$, $m_2=5$, and in the last $4$
cases, we set $t_1=\pm 1$, $t_2=\pm 2$, $t_3=\pm 3$, $t_4=\pm 5$,
and $m_1=7$. In this way, we obtain the connected component for
each $D_{\sigma}^-(h)$. To denote a choice of sign we write a
list of $6$ symbols, with $+$ or $-$ indicating the sign
associated to an element $j$ of $I(\sigma)$, $0$ indicating an
element of $J(\sigma)$, and $*$ indicating elements of $K(\sigma)$.

\begin{center}
1. $y_{1x12x2}(t_1,m_1,t_2,m_2)$.

\vspace*{0.3cm}
\begin{tabular}{c|c|c}
Sign choice & Signs of $a,b,c,d,e,f$ & Connected Component \\
\hline
$0+*0+*$ & $---+++$ & $K\longleftrightarrow 11$ \\
$0+*0-*$ & $---+-+$ & $J\longleftrightarrow 9$ \\
$0-*0+*$ & $--+-++$ & $I\longleftrightarrow 10$ \\
$0-*0-*$ & $--+--+$ & $K\longleftrightarrow 11$
\end{tabular}
\vspace*{0.3cm}
\end{center}

\begin{center}
2. $y_{xx1x1x}(t_1,t_2,t_3,m_1, t_4)$.

\vspace*{0.3cm}
\begin{tabular}{c|c|c}
Sign choice & Signs of $a,b,c,d,e,f$ & Connected Component \\
\hline
$++0+*+$ & $+-+++-$ & $H\longleftrightarrow 8 $\\
$++0+*-$ & $--+++-$ & $K\longleftrightarrow 11$ \\
$++0-*+$ & $+-++++$ & $I\longleftrightarrow 10$ \\
$++0-*-$ & $--++++$ & $F\longleftrightarrow 6 $\\
$+-0+*+$ & $+--+-+$ & $F\longleftrightarrow 6 $\\
$+-0+*-$ & $---+-+$ & $J\longleftrightarrow 9 $\\
$+-0-*+$ & $+--+--$ & $K\longleftrightarrow 11$ \\
$+-0-*-$ & $---+--$ & $G\longleftrightarrow 7 $\\
$-+0+*+$ & $+-+---$ & $J\longleftrightarrow 9 $\\
$-+0+*-$ & $--+---$ & $E\longleftrightarrow 5 $\\
$-+0-*+$ & $+-+--+$ & $H\longleftrightarrow 8 $\\
$-+0-*-$ & $--+--+$ & $K\longleftrightarrow 11$ \\
$--0+*+$ & $+---++$ & $K\longleftrightarrow 11$ \\
$--0+*-$ & $----++$ & $G\longleftrightarrow 7 $\\
$--0-*+$ & $+---+-$ & $E\longleftrightarrow 5 $\\
$--0-*-$ & $----+-$ & $I\longleftrightarrow 10$
\end{tabular}
\vspace*{0.3cm}
\end{center}

\begin{center}
3. $y_{1x1xxx}(t_1,m_1,t_2,t_3,t_4)$.

\vspace*{0.3cm}
\begin{tabular}{c|c|c}
Sign choice & Signs of $a,b,c,d,e,f$ & Connected Component \\
\hline
$0+*+++$ & $+--+-+$ & $F\longleftrightarrow 6 $\\
$0+*++-$ & $---+-+$ & $J\longleftrightarrow 9 $\\
$0+*+-+$ & $+-+--+$ & $H\longleftrightarrow 8 $\\
$0+*+--$ & $--+--+$ & $K\longleftrightarrow 11$ \\
$0+*-++$ & $+-++-+$ & $K\longleftrightarrow 11$ \\
$0+*-+-$ & $--++-+$ & $G\longleftrightarrow 7 $\\
$0+*--+$ & $+----+$ & $J\longleftrightarrow 9 $\\
$0+*---$ & $-----+$ & $E\longleftrightarrow 5 $\\
$0-*+++$ & $+--+++$ & $H\longleftrightarrow 8 $\\
$0-*++-$ & $---+++$ & $K\longleftrightarrow 11$ \\
$0-*+-+$ & $+-+-++$ & $E\longleftrightarrow 5 $\\
$0-*+--$ & $--+-++$ & $I\longleftrightarrow 10$ \\
$0-*-++$ & $+-++++$ & $I\longleftrightarrow 10$ \\
$0-*-+-$ & $--++++$ & $F\longleftrightarrow 6 $\\
$0-*--+$ & $+---++$ & $K\longleftrightarrow 11$ \\
$0-*---$ & $----++$ & $G\longleftrightarrow 7$
\end{tabular}
\vspace*{0.3cm}
\end{center}

\begin{center}
4. $y_{12x21x}(t_1,m_1,m_2, t_2)$.

\vspace*{0.3cm}
\begin{tabular}{c|c|c}
Sign choice & Signs of $a,b,c,d,e,f$ & Connected Component \\
\hline
$00+**+$ & $+-+++-$ & $H\longleftrightarrow 8 $\\
$00+**-$ & $--+++-$ & $K\longleftrightarrow 11$ \\
$00-**+$ & $+--+--$ & $K\longleftrightarrow 11$ \\
$00-**-$ & $---+--$ & $G\longleftrightarrow 7$
\end{tabular}
\vspace*{0.3cm}
\end{center}

\begin{center}
5. $y_{xxx2x2}(t_1,t_2,t_3,t_4,m_1)$.

\vspace*{0.3cm}
\begin{tabular}{c|c|c}
Sign choice & Signs of $a,b,c,d,e,f$ & Connected Component \\
\hline
$+++0+*$ & $-+++-+$ & $F\longleftrightarrow 6$ \\
$+++0-*$ & $-+++++$ & $H\longleftrightarrow 8$ \\
$++-0+*$ & $--++--$ & $I\longleftrightarrow 10$ \\
$++-0-*$ & $--+++-$ & $K\longleftrightarrow 11$ \\
$+-+0+*$ & $-+---+$ & $G\longleftrightarrow 7 $\\
$+-+0-*$ & $-+-+--$ & $F\longleftrightarrow 6 $\\
$+--0+*$ & $---+++$ & $K\longleftrightarrow 11$ \\
$+--0-*$ & $---+-+$ & $J\longleftrightarrow 9 $\\
$-++0+*$ & $-++-++$ & $K\longleftrightarrow 11$ \\
$-++0-*$ & $-++--+$ & $J\longleftrightarrow 9 $\\
$-+-0+*$ & $--+-+-$ & $H\longleftrightarrow 8 $\\
$-+-0-*$ & $--+---$ & $E\longleftrightarrow 5 $\\
$--+0+*$ & $-+----$ & $I\longleftrightarrow 10$ \\
$--+0-*$ & $-+--+-$ & $K\longleftrightarrow 11$ \\
$---0+*$ & $-----+$ & $E\longleftrightarrow 5 $\\
$---0-*$ & $---+--$ & $G\longleftrightarrow 7$
\end{tabular}
\vspace*{0.3cm}
\end{center}

\begin{center}
6. $y_{x2x2xx}(t_1,t_2,m_1,t_3,t_4)$.

\vspace*{0.3cm}
\begin{tabular}{c|c|c}
Sign choice & Signs of $a,b,c,d,e,f$ & Connected Component \\
\hline
$+0+*++$ & $--++--$ & $I\longleftrightarrow 10$ \\
$+0+*+-$ & $-+-+--$ & $F\longleftrightarrow 6 $\\
$+0+*-+$ & $-+++--$ & $K\longleftrightarrow 11$ \\
$+0+*--$ & $---+--$ & $G\longleftrightarrow 7 $\\
$+0-*++$ & $--++++$ & $F\longleftrightarrow 6 $\\
$+0-*+-$ & $-+-+++$ & $J\longleftrightarrow 9 $\\
$+0-*-+$ & $-+++++$ & $H\longleftrightarrow 8 $\\
$+0-*--$ & $---+++$ & $K\longleftrightarrow 11$ \\
$-0+*++$ & $--+-+-$ & $H\longleftrightarrow 8 $\\
$-0+*+-$ & $-+--+-$ & $K\longleftrightarrow 11$ \\
$-0+*-+$ & $-++-+-$ & $E\longleftrightarrow 5 $\\
$-0+*--$ & $----+-$ & $I\longleftrightarrow 10$ \\
$-0-*++$ & $--+--+$ & $K\longleftrightarrow 11$ \\
$-0-*+-$ & $-+---+$ & $G\longleftrightarrow 7 $\\
$-0-*-+$ & $-++--+$ & $J\longleftrightarrow 9 $\\
$-0-*--$ & $-----+$ & $E\longleftrightarrow 5$
\end{tabular}
\vspace*{0.3cm}
\end{center}

\begin{center}
7. $y_{x21x12}(t_1,t_2,m_1,m_2)$.

\vspace*{0.3cm}
\begin{tabular}{c|c|c}
Sign choice & Signs of $a,b,c,d,e,f$ & Connected Component \\
\hline
$+00+**$ & $-+++-+$ & $F\longleftrightarrow 6$ \\
$+00-**$ & $-+++--$ & $K\longleftrightarrow 11$ \\
$-00+**$ & $---+++$ & $K\longleftrightarrow 11$ \\
$-00-**$ & $---++-$ & $E\longleftrightarrow 5$
\end{tabular}
\vspace*{0.3cm}
\end{center}

We now have a decomposition of each connected component of
$\mathcal B^*$ into a disjoint union of subsets of form
$D_{\sigma}(h)$. This allows us to calculate the Euler
characteristic (for compactly supported cohomology) of each
connected component $X$ as an alternating sum:

 $$
\chi(X)=\sum_{c=0}^{\dim(X)} (-1)^c n_c,
 $$
where $n_c$ is the number of subsets of form $D_{\sigma}(h)$ of
codimension $c$ contained in $X$ (note that the dimension of $X$ is even).
Since the components $X$ are
smooth (open in $\mathcal B$) this compactly supported Euler
characteristic coincides with the usual one by Poincar\'e duality.
We note that each $D_{\sigma}(h)$ has codimension at most $2$.

\begin{center}
\begin{figure}
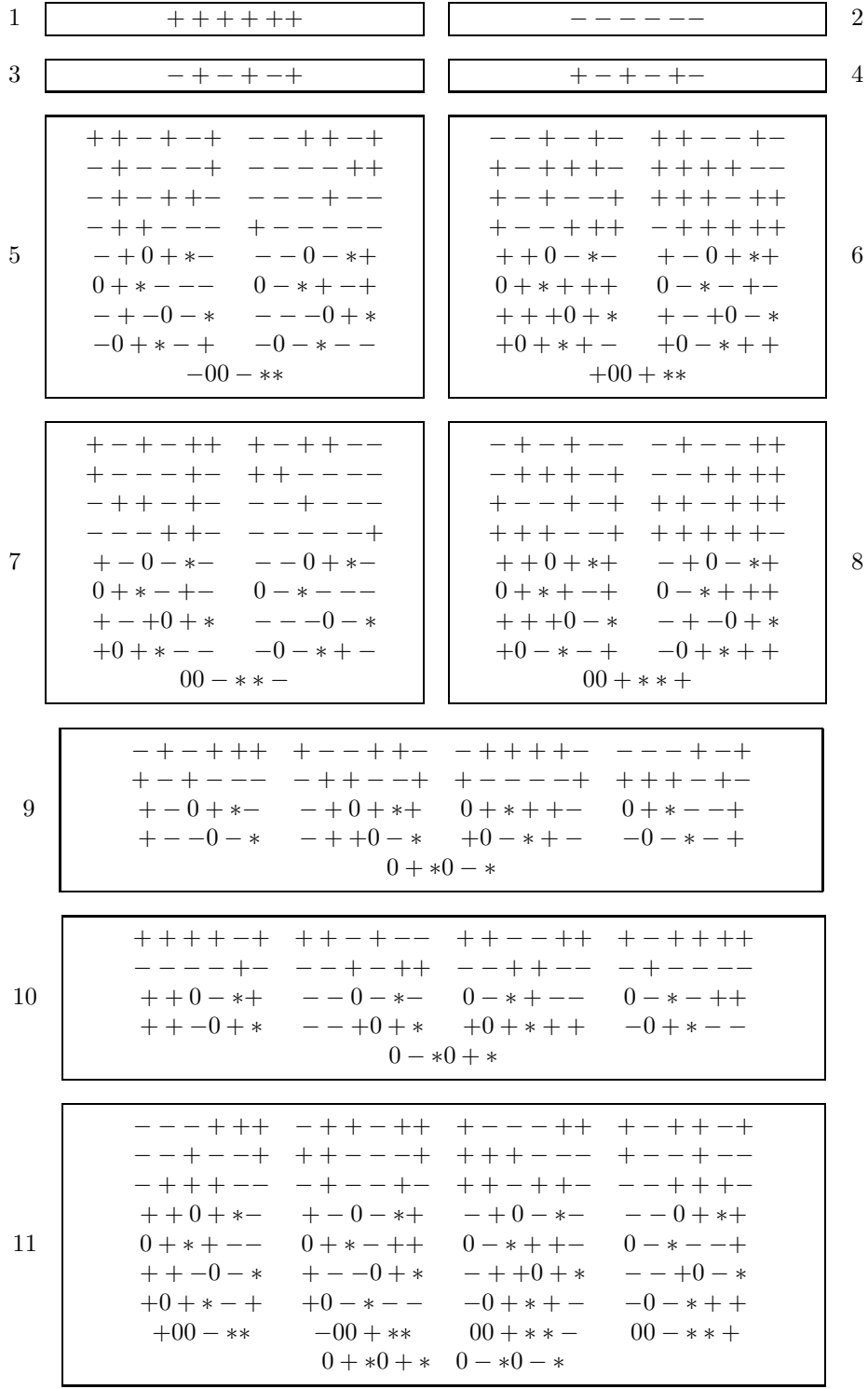

\begin{tabular}{cccc}
\ 1 &
\framebox[5.4cm]{$++++++$}
&
 \framebox[5.4cm]{$------$}
& 2
\end{tabular}

\vspace*{0.3cm}

\begin{tabular}{cccc}
\ 3 &
 \framebox[5.4cm]{$-+-+-+$}
&
 \framebox[5.4cm]{$+-+-+-$}
& 4
\end{tabular}

\vspace*{0.3cm}

\begin{tabular}{cccc}
\ 5 &
 \framebox[5.4cm]{
$\begin{array}{c}
\begin{array}{cc}
++-+-+ & --++-+ \\
-+---+ & ----++ \\
-+-++- & ---+-- \\
-++--- & +----- \\
-+0+*- & --0-*+ \\
0+*--- & 0-*+-+ \\
-+-0-* & ---0+* \\
-0+*-+ & -0-*--
\end{array}\\
-00-**
\end{array}$
}
&
 \framebox[5.4cm]{
$\begin{array}{c}
\begin{array}{cc}
--+-+- & ++--+- \\
+-+++- & ++++-- \\
+-+--+ & +++-++ \\
+--+++ & -+++++ \\
++0-*- & +-0+*+ \\
0+*+++ & 0-*-+- \\
+++0+* & +-+0-* \\
+0+*+- & +0-*++
\end{array}\\
+00+**
\end{array}$
}
& 6
\end{tabular}

\vspace*{0.3cm}

\begin{tabular}{cccc}
\ 7 &
 \framebox[5.4cm]{
$\begin{array}{c}
\begin{array}{cc}
+-+-++ & +-++-- \\
+---+- & ++---- \\
-++-+- & --+--- \\
---++- & -----+ \\
+-0-*- & --0+*- \\
0+*-+- & 0-*--- \\
+-+0+* & ---0-* \\
+0+*-- & -0-*+-
\end{array}\\
00-**-
\end{array}$
}
&
 \framebox[5.4cm]{
$\begin{array}{c}
\begin{array}{cc}
-+-+-- & -+--++ \\
-+++-+ & --++++ \\
+--+-+ & ++-+++ \\
+++--+ & +++++- \\
++0+*+ & -+0-*+ \\
0+*+-+ & 0-*+++ \\
+++0-* & -+-0+* \\
+0-*-+ & -0+*++
\end{array}\\
00+**+
\end{array}$
}
& 8
\end{tabular}

\vspace*{0.3cm}

\begin{tabular}{ccc}
\ 9 &
\framebox[10.9cm]{
$\begin{array}{c}
\begin{array}{cccc}
-+-+++ & +--++- & -++++- & ---+-+ \\
+-+--- & -++--+ & +----+ & +++-+- \\
+-0+*- & -+0+*+ & 0+*++- & 0+*--+ \\
+--0-* & -++0-* & +0-*+- & -0-*-+
\end{array}\\
0+*0-*
\end{array}$
}
& \  \end{tabular}

\vspace*{0.3cm}

\begin{tabular}{ccc}
10 &
\framebox[10.9cm]{
$\begin{array}{c}
\begin{array}{cccc}
++++-+ & ++-+-- & ++--++ & +-++++ \\
----+- & --+-++ & --++-- & -+---- \\
++0-*+ & --0-*- & 0-*+-- & 0-*-++ \\
++-0+* & --+0+* & +0+*++ & -0+*--
\end{array}\\
0-*0+*
\end{array}$
}
& \
\end{tabular}

\vspace*{0.3cm}

\begin{tabular}{ccc}
11 &
\framebox[10.9cm]{
$\begin{array}{c}
\begin{array}{cccc}
---+++ & -++-++ & +---++ & +-++-+ \\
--+--+ & ++---+ & +++--- & +--+-- \\
-+++-- & -+--+- & ++-++- & --+++- \\
++0+*- & +-0-*+ & -+0-*- & --0+*+ \\
0+*+-- & 0+*-++ & 0-*++- & 0-*--+ \\
++-0-* & +--0+* & -++0+* & --+0-* \\
+0+*-+ & +0-*-- & -0+*+- & -0-*++ \\
+00-** & -00+** & 00+**- & 00-**+ \\
\end{array}\\
\begin{array}{cc}
0+*0+* & 0-*0-*
\end{array}
\end{array}$
} & \  \end{tabular} \caption{The connected components of
$\mathcal{B}^*$ with their cell decompositions.\label{components}}
\end{figure}
\end{center}

All the required information is displayed in
Figure~\ref{components}, using a format similar to that used
in~\cite{Rie:BruCel}, grouping the subsets of form
$D_{\sigma}(h)$ by connected component. We denote each
$D_{\sigma}(h)$ by the sign choice string listed above (first
column); note that the number of $0$'s gives the codimension.

\begin{rems}
Fix $\sigma$ such that $D_{\sigma}$ is of codimension $2$ (cases
$1$--$3$). Then, as $|I(\sigma)|=2$ there are $2^2=4$ possible
choices of sign, giving rise to $4$ subsets $D_{\sigma}(h)$. Two
of these lie in component $11$, and the other two lie in either
components $5$ and $6$, components $7$ and $8$, or components $9$
and $10$ respectively.

Fix $\sigma$ such that $D_{\sigma}$ is of codimension $1$ (cases
$4$--$7$). Then, as $|I(\sigma)|=4$ there are $2^4=16$ possible
choices of sign, giving rise to $16$ subsets $D_{\sigma}(h)$.
Each connected component $5$--$10$ contains precisely two of
these subsets, while connected component $11$ contains precisely
$4$ of these.
\end{rems}

\begin{thm}
The Euler characteristic of each connected component of $\mathcal{B}^*$
%U^-\cap B^+\dot w_0 B^+$
is given in Table II. We also give, for each connected component,
and for codimension $0$, $1$, and $2$, the number of subsets of
the form $D_{\sigma}(h)$ of that codimension contained in that
component. We thus conclude that there are $10$ connected
components of Euler characteristic $1$ and one of Euler
characteristic $2$, confirming a total Euler characteristic of
$\mathcal{B}^*$ of $12$, the sum of these.
\end{thm}

\vspace*{0.3cm}

\begin{center}
Table II. Euler characteristics of the connected components.

\vspace*{0.3cm}

\begin{tabular}{c|c|c|c|c}
Connected component & Codim 0 & Codim 1 & Codim 2 & Euler characteristic \\
\hline
 1 & 1 & 0 & 0 & 1 \\
 2 & 1 & 0 & 0 & 1 \\
 3 & 1 & 0 & 0 & 1 \\
 4 & 1 & 0 & 0 & 1 \\
 5 & 8 & 8 & 1 & 1 \\
 6 & 8 & 8 & 1 & 1 \\
 7 & 8 & 8 & 1 & 1 \\
 8 & 8 & 8 & 1 & 1 \\
 9 & 8 & 8 & 1 & 1 \\
10 & 8 & 8 & 1 & 1 \\
11 & 12 & 16 & 6 & 2
\end{tabular}
\end{center}

% ----------------------------------------------------------------
\bibliographystyle{amsplain}
\bibliography{biblio}
\end{document}